\documentclass[11pt]{article}

\usepackage{tikz}
\usetikzlibrary{decorations.pathreplacing}
\usepackage{epsfig }
\usepackage{latexsym}
\usepackage{amsfonts,amsmath,amssymb,textcomp}
\allowdisplaybreaks

\usepackage{amsfonts,amsmath,amssymb}

\usepackage{epsfig}
\usepackage{latexsym}

\usepackage{algorithm}
\usepackage{algorithmic}

\newtheorem{thm}{Theorem}[section]
 \newtheorem{lem}[thm]{Lemma}
\newtheorem{cor}[thm]{Corollary}
\newtheorem{prop}[thm]{Proposition}
\newtheorem{rem}[thm]{Remark}
\newtheorem{deff}[thm]{Definition}
\newtheorem{coni}[thm]{Coniecture}
\newcommand{\bconi}{\begin{coni}}
\newcommand{\econi}{\end{coni}}
\newcommand{\bth}{\begin{thm}}
\newcommand{\ethGL}{\end{thm}}
\newcommand{\bl}{\begin{lem}}
\newcommand{\el}{\end{lem}}
\newcommand{\bdf}{\begin{deff}}
\newcommand{\edf}{\end{deff}}
\newcommand{\bcor}{\begin{cor}}
\newcommand{\ecor}{\end{cor}}
\newcommand{\bprop}{\begin{prop}}
\newcommand{\eprop}{\end{prop}}
\newcommand{\brem}{\begin{rem}}
\newcommand{\erem}{\end{rem}}
\newcommand{\beq}{\begin{equation}}
\newcommand{\eeq}{\end{equation}}
\newcommand{\beqn}{\begin{eqnarray}}
\newcommand{\eeqn}{\end{eqnarray}}
\newcommand{\beqns}{\begin{eqnarray*}}
\newcommand{\eeqns}{\end{eqnarray*}}

\newcommand{\BO}{\mathcal{O}}
\newcommand{\BN}{\mathcal{N}}

\newcommand{\non}{\nonumber}
\newcommand{\ba}{\begin{array}}
\newcommand{\ea}{\end{array}}
\newcommand{\bit}{\begin{itemize}}
\newcommand{\eit}{\end{itemize}}
\newcommand{\ben}{\begin{enumerate}}
\newcommand{\een}{\end{enumerate}}
\newcommand{\bal}{\begin{align}}
\newcommand{\bals}{\begin{align*}}
\newcommand{\bs}{\begin{skip}}
\newcommand{\eal}{\end{align}}
\newcommand{\eals}{\end{align*}}
\newcommand{\es}{\end{skip}}
\newcommand{\babs}{\begin{abstract}}
\newcommand{\eabs}{\end{abstract}}

\newcommand{\B}{\quad}

\newcommand{\kl}{[\hspace{-0.5 mm}[}
\newcommand{\kr}{]\hspace{-0.5 mm}]}
\newcommand{\lb}{\left [}
\newcommand{\rb}{\right ]}
\newcommand{\lp}{\left (}
\newcommand{\rp}{\right )}

\newcommand{\lf}{\lfloor}
\newcommand{\rf}{\rfloor}
\newcommand{\dq}{\doteq}
\newcommand{\QB}[2]{\genfrac{[}{]}{0pt}{}{#1}{#2}_q}

\def\ii{\mathbf{i}}

\newcommand{\srf}{\sqrt{5}}

\newcommand{\Rb}{\overline{R}}
\newcommand{\Eb}{\overline{\E}}
\newcommand{\mub}{\overline{\mu}}
\newcommand{\sigb}{\overline{\sig}}

\newcommand{\At}{\tilde{A}}
\newcommand{\Bt}{\tilde{B}}

\newcommand{\de}{\delta}
\newcommand{\De}{\Delta}

\newcommand{\la}{\lambda}

\newcommand{\al}{\alpha}
\newcommand{\be}{\beta}
\newcommand{\ra}{\rightarrow}

\newcommand{\srn}{\sqrt{n}}
\newcommand{\srm}{\sqrt{m}}

\newcommand{\Fi}{\varphi}
\newcommand{\FI}{\phi}
\newcommand{\eps}{\varepsilon}
\newcommand{\II}{\infty}
\newcommand{\tet}{\theta}
\newcommand{\sig}{\sigma}
\newcommand{\sigd}{\sigma^2}

\newcommand{\gam}{\gamma}
\newcommand{\Gam}{\Gamma}
\newcommand{\ta}{\tau}
\newcommand{\ro}{\rho}

\newcommand{\conv}{\stackrel{\mathcal{D}}{\sim}}
\newcommand{\ds}{\stackrel{ .}{\sim}}

\newcommand{\E}{\mbox{$\mathbb E$}}
\newcommand{\V}{\mbox{$\mathbb V$}}
\def\P{{\mathbb {P}}}

\newcommand{\bin}[2]
{
{#1\choose #2}
}

\def\la{\lambda}

\def\St2#1#2{\left\{\ba{c}#1\\#2\ea\right\}}

\title{\bf Some large polyominoes' perimeter: a stochastic analysis }
\author{Guy Louchard\thanks{Universit\'e Libre de Bruxelles,
D\'epartement d'Informatique, CP 212, Boulevard du Triomphe, B-1050
Bruxelles, Belgium, email: louchard@ulb.ac.be}
}

\date{\today}
\begin{document}
\maketitle
\babs
In this paper, we analyze the stochastic properties of some large size (area) polyominoes' perimeter such that  
the directed column-convex polyomino, the column-convex polyomino, the directed diagonally-convex polyomino, the staircase (or parallelogram) polyomino, the escalier polyomino,  the wall (or bargraph) polyomino.
All polyominoes considered here are made of contiguous, not-empty columns, without holes, such that each column  must be adjacent to some cell of the previous column. We compute the asymptotic (for large size $n$)  Gaussian distribution of the perimeter, including the corresponding Markov property of the chain of columns, and the convergence to classical Brownian motions of the perimeter seen as a trajectory according to the successive columns. All polyominoes of size $n$ are considered as equiprobable. 
\eabs


\section{Introduction}
In this paper, a polyomino is a set of cells on a square lattice such that every cell of the polyomino can be reached from any other cell by a sequence of cells of the polyomino. The perimeter ( the length of the border) has been the subject of a large literature. We will not provide all references, we refer to the rather complete bibliography given in  Bousquet-M\'elou \cite{BM93},\cite{BM96},
Bousquet-M\'elou  and Brak \cite{BM08}. Let us also add  Fereti\'c and Svrtan \cite{FE96}, Fereti\'c \cite{FE02},
 Delest and F\'edou \cite{DF89}, Blecher et al. \cite{BBKM17}, Louchard \cite{GL18}. 

In this paper, we analyze the stochastic properties of some large size (area) polyominoes' perimeter such that  (precise definitions are given in the text)
the directed column-convex polyomino (dcc), the column-convex polyomino (cc), the directed diagonally-convex polyomino (dc), the staircase (or parallelogram) polyomino (st), the escalier polyomino (es),  the wall (or bargraph) polyomino (wa).
All polyominoes considered here are made of contiguous, not-empty columns, without holes, such that each column (of size $j$) must be adjacent to some cell of the previous column (of size $k$). We will denote by $U(k,j)$ (characterizing each polyomino) the possibility function giving the number of ways of gluing the two
columns together. We compute the asymptotic (for large size $n$) Gaussian  distribution of the perimeter, including the corresponding Markov property of the chain of columns, and the convergence to classical Brownian motions of the perimeter seen as a trajectory according to the successive columns. 
This confirms   the ``filament silhouette" of the structures, that had been observed by previous simulations. As said in Flajolet and Sedgewick \cite [p.662]{FLSe09} ,``a random parallelogram is most likely to resemble a slanted stack of fairly short segments". This is proved here for all our polyominoes. All polyominoes of size $n$ are considered as equiprobable. The first five ones are treated with similar methods, the bargraph is analyzed with a different technique.

\section{The dcc perimeter}
A directed column convex polyomino (dcc) is made of contiguous columns such that the base cell of each column must be adjacent to some cell of the previous column. 
We have partially considered this polyomino in \cite{GL96}. In this section, besides the perimeter's analysis,  we refine the polyomino
stochastic description with another technique, that we will use in all other polyominoes. So we explain it in great detail in this section and provide all necessary notations and computations we use in the sequel.

 For dcc, the gluing  function is given by $U(k,j)=k$. Note that it depends only on $k$, this will not be the case for our following 
polyominoes.
In this paper, we denote by $\FI(w,\tet,z)$ the three-dimension  generating function (GF) where $z$ marks the polyomino's size $n$ (area), $w$ marks the number $m$ of columns (width)  and $\tet$ marks the size $j$  of the last column. All other interesting  parameters and stochastic distributions are related to $\FI(w,\tet,z)$. 
In the following subsections, we will first consider $\FI(w,\tet,z)$ and its derived properties, then the Markov chain corresponding to the dcc, next the perimeter conditioned on the number of columns $m$ and finally the perimeter conditioned on the size $n$. Asymptotic relations always means when $n\ra \II$.
\subsection{The generating functions}                                        \label{S21}
To compute the generating functions we will use the ``adding of a slice" technique, initiated by Temperley \cite {TE56}, popularized   by many combinatorists (see for instance Bousquet-M\'elou \cite{BM96}) and summarized in
 Flajolet and Sedgewick \cite[p.366]{FLSe09}. The analysis we apply here to dcc has already been initialized in \cite{GL297} for cc, but for the sake of completeness, we present it again  here, with some complements. 
 
 We denote by $T(m,n,j,\ell)$ the total number of dcc with area $n$, width $m$, last column size $j$ and first column size $\ell$ (similar notations for partially parametrized $T()$.

For any function $g(\tet,..)$, set:
\[g'(\tet,..)\dq \partial_\tet g(\tet,..).\]
When we use the symbol $\dq$ (and similarly for $\ds$), this corresponds to a relation valid for \emph{all} polyominoes, otherwize, the usual $:=$ corresponds to the particular polyomino under consideration.
Denote by $\Fi(m,j,z)$ the GF corresponding to polyominoes with $m$ columns, last column of size $j$ and any first column's size. We have
\[\Fi(m,j,z)=\sum_{k=1} ^\II k \Fi(m-1,k,z).\]
Now we mark $j$ by $\tet$. This leads to 
\bals
 \Gam(m,\tet,z)&\dq \sum_{j=1} ^\II \tet^j \Fi(m,j,z)= \sum_{j=1} ^\II \tet^j z^j \sum_{k=1} ^\II k \Fi(m-1,k,z)
 = \tilde{f}_2(\tet,z)\Gam'(m-1,1,z),\B \tilde{f}_2(\tet,z)= \tet z/(1-\tet z),\\
 \Gam(1,\tet,z)&\dq \tilde{f}_0(\tet,z),\tilde{f}_0(\tet,z)= \tilde{f}_2(\tet,z).
\end{align*}
Set
\bals
\Gam(m,\tet,z)&\dq z^m \tet  \De(m,\tet,z),\Gam'(m,\tet,z)\dq z^m   \De(m,\tet,z)+z^m \tet  \De'(m,\tet,z),\mbox{ hence }\\
   \De(m,\tet,z)&\dq f_1(\tet,z)\De(m-1,1,z)+f_2(\tet,z)\De'(m-1,1,z), m\geq 2, \B f_2(\tet,z)= 1/(1-\tet z),f_1(\tet,z)=  f_2(\tet,z),
	\\
 \De(1,\tet,z)&\dq f_0(\tet,z), f_0(\tet,z)=  f_2(\tet,z). 
\end{align*}
Define
\[\psi(\xi,\tet,z)\dq \sum_{m=1} ^\II \xi^m \De(m,\tet,z).\]
We obtain
\bals
  \psi(\xi,\tet,z)&\dq \xi[f_0(\tet,z)+f_1(\tet,z) D_1(\xi,z)+f_2(\tet,z) D_2(\xi,z)],\mbox{ with }\\
  D_1(\xi,z)&\dq \psi(\xi,1,z),D_2(\xi,z)\dq \psi'(\xi,1,z),\mbox{ and }\\
	\FI(w,\tet,z)&\dq\sum_{w=1}^\II w^m \Gam(m,\tet,z)\dq  \tet\psi(wz,\tet,z).
\end{align*}	
To obtain $D_1,D_2$, we compute
\bal
\psi'(\xi,\tet,z)&\dq\xi[f'_0(\tet,z)+f'_1(\tet,z) D_1(\xi,z)+f'_2(\tet,z) D_2(\xi,z)],\non\\ 
	D_1(\xi,z)&\dq\xi[f_0(1,z)+f_1(1,z) D_1(\xi,z)+f_2(1,z) D_2(\xi,z)],                      \label{E1}\\
	D_2(\xi,z)&\dq\xi[f'_0(1,z)+f'_1(1,z) D_1(\xi,z)+f'_2(1,z) D_2(\xi,z)].              \label{E2}       
\end{align}
Solving, we get

\bals
	D_1(\xi,z)&\dq \frac{N_1(\xi,z)}{h(\xi,z)},D_2(\xi,z)\dq \frac{N_2(\xi,z)}{h(\xi,z)},\\
	N_1(\xi,z)&= \xi(z-1),N_2(\xi,z)= -\xi z,h(\xi,z)= -z^2+2z+\xi-1,\mbox{ and setting } \xi=zw,\mbox{ this gives }\\
	N_1(w,z)&= zw(z-1),N_2(w,z) -z^2w,h(w,z) -z^2+2z+zw-1.
	\end{align*}
	\[\mbox{ The root of smallest module of } h(1,\ro)=0 \mbox{ is given by } \ro= \frac{3-\srf}{2}= \frac{1}{\FI^2},
	(\FI \mbox{ is the golden ratio }).\]
	\[\FI(w,\tet,z)\dq \tet\psi(wz,\tet,z)\dq w\tet z[f_0(\tet,z)+f_1(\tet,z) D_1(w,z)+f_2(\tet,z) D_2(w,z)],\]
		
		$f_0(\tet,z)$ corresponds of course to the first column. In order not to burden the notations, we use indifferently
		$F(\xi),F(w)$, where $\xi=zw$, clearly depending on the context.  Note that we recover $\ro$, already computed in \cite{GL96}, in a simple way.
		
		If we set $\tet=1$ in $\FI$, we have $D_1(w,z)$. When $w=1$, this gives the GF of the total number $T(.,n)$ of size $n$ dcc. By classical singularity singularity analysis, this leads to
		\[T(.,n)\ds \frac{C_1}{\ro^n},n\ra\II,C_1\dq -\frac{N_1(1,\ro)}{\ro h_z(1,\ro)}= \frac12-\frac{\srf}{10}.\]
		But we get more. By Bender's theorems $1$ and $3$ in \cite {BE73} \footnote{See Appendix \ref{A4}}, (see also Flajolet and Sedgewick, \cite [Thm.IX.9]{FLSe09}) we derive the asymptotic distribution of the width $M$, given  area $n$. Let  $F_z$ mean differentiation of $F$ w.r.t. $z$, and similarly for other notations. Set
		\bals
		r_1&\dq -h_w/h_z= -\frac{z}{-2z+2+w},\\
	r_2&\dq -(r_1^2h_{zz}+2r_1h_{zw}+h_w+h_{ww})/h_z= -z(-6z+2w+4z^2-4zw+w^2)/(-2z+2+w)^3,\\
	 \mu_1&\dq -r_1/\ro= \frac{\srf}{5},\sig_1^2\dq \mu_1^2-r_2/\ro= 2\frac{\srf}{25}.
	\end{align*}
	Then Bender's theorems lead to 
		\bth
		The width $M$ of a dcc of large given area $n$ is asymptotically Gaussian:
		\footnote{$\conv$ means convergence in distribution}
		\beq \frac{M-n\mu_1}{\srn\sig_1}\conv \mathcal {N}(0,1)=\ta_0,\mbox{ say },n\ra\II,                   \label{E3}\eeq
		also a local limit theorem holds:
		\[T(m,n)\ds \frac{C_1}{\ro^n}\frac{e^{-(m-n\mu_1)^2/(2n\sig_1^2)}}{\sqrt{2\pi n}\sig_1},n\ra\II,m-n\mu_1=\BO(\srn).\]
		\ethGL
		The verification of condition (V) of Bender's Theorem $3$ (which is essential to go from a central limit theorem to a local limit theorem) is easy: the function $h(w,z)$ has the following property: $h(e^s,z)$ is analytic and bounded for
		\[|z|\leq |r(\Re(s))|(1+\de)\mbox{ and }\eps\leq|\Im(s)|\leq \pi\]
		for some $\eps>0,\de>0$, where $r(s)$ is the suitable solution of the equation $h(e^s,r(s))=0$ (i.e. with $r(0)=\ro$). This will be valid for all functions $h(w,z)$ used in the following sections.
		
	Now if we fix $m$ and consider $n$ as a variable (there are, of course, an infinite number of dcc for a given $m$), we can obtain another asymptotic expression for $T(m,n)$. The conditioned distribution is given by 
		\[[ w^m z^n]D_1(w,z)\dq \frac{1}{\ro^n}[ w^m z^n]D_1(w,\ro z).\]
					For $z=1$, the dominant singularity of $D_1(w,\ro )$ is $w=1$. With 
					\[C_2\dq -\frac{N_1(1,\ro)}{ h_w(1,\ro)}= \frac{\srf-1}{2} ,
	C_2\dq C_1/\mu_1,\mu_2\dq 1/\mu_1=\srf,\sig_2^2\dq \sig_1^2/\mu_1^3= 2.\]
	We derive the following theorem
	\bth                                                                          \label{T22}
	For a large given width $m$, $T(m,n)$ is asymptotically given by
	\[T(m,n)\ds \frac{C_2}{\ro^n}\frac{e^{-(n-m\mu_2)^2/(2m\sig_2^2)}}{\sqrt{2\pi m}\sig_2},m\ra\II,n-n\mu_2=\BO(\srm).\]
	\ethGL
Let us now turn to the asymptotic GF $G(\tet)$ of the last column size.	We first compute
\[[z^n]\FI(1,\tet,z)\ds \frac{\tet \ro [f_1(\tet,\ro) N_1(1,\ro)+f_2(\tet,\ro) N_2(1,\ro)]}{-\ro^n\ro h_z(1,\ro)},n\ra \II\]
uniformly for $\tet$ in some complex neighbourhood of the origin. This may be checked by the method of singularity analysis of Flajolet and Odlyzko, as used in Flajolet and Soria \cite{FS93}.

 Normalizing by $T(.,n)$, this leads to the following theorem
\bth
\bals
G(\tet)&\dq \tet\ro\lb f_1(\tet,\ro)+f_2(\tet,\ro)\frac{N_2(1,\ro)}{N_1(1,\ro)}\rb,\\
	\pi(j)&\dq [\tet^j]G(\tet)= \ro^j\lp 1+\frac{N_2(1,\ro)}{N_1(1,\ro)}\rp = \ro^j\frac{\srf+1}{2}.
\end{align*}
\ethGL
But multiplying (\ref{E1}),(\ref{E2}) by $h(w,z)$ and letting $w\ra 1,z\ra \ro$, we have
\bal
N_1(1,\ro)&\dq \ro[f_1(1,\ro) N_1(1,\ro)+f_2(1,\ro) N_2(1,\ro)],                            \label{E4}\\
	N_2(1,\ro)&\dq \ro[f'_1(1,\ro) N_1(1,\ro)+f'_2(1,\ro) N_2(1,\ro)].                           \non
\end{align}
Hence
\[ G(1)\dq 1,G'(1)\dq 1+\frac{N_2(1,\ro)}{N_1(1,\ro)}= \frac{\srf+1}{2}.\]
Let us analyze the asymptotic distribution of the last column size in a  dcc of large area $n$ and width $m$. Again Bender's theorems lead to the following GF (the notation is clear here)
\bals
 T(m,n,\tet)&\ds\frac{e^{-(m-n\mu_1)^2/(2n\sig_1^2)}}{\sqrt{2\pi n}\sig_1}\frac{1}{-\ro^n\ro h_z(1,\ro)}\tet \ro [f_1(\tet,\ro) N_1(1,\ro)+f_2(\tet,\ro) N_2(1,\ro)],\mbox{ or }\\
 T(m,n,\tet)&\ds \frac{C_1}{\ro^n}\frac{e^{-(m-n\mu_1)^2/(2n\sig_1^2)}}{\sqrt{2\pi n}\sig_1}G(\tet)\dq \frac{C_2}{\ro^n}\frac{e^{-(n-m\mu_2)^2/(2m\sig_2^2)}}{\sqrt{2\pi m}\sig_2}G(\tet),\\
n&\ra\II,m-n\mu_1=\BO(\srn),m\ra\II,n-n\mu_2=\BO(\srm).
 \end{align*}
 We now turn to the case where the first column possesses $i$ cells. This leads to (note that $h(w,z)$ remains the same, independently of $i$)
 
	\bals
 \Gam(1,\tet,z)&\dq \tet^iz^i,f_0(\tet,z,i)\dq \tet^{i-1}z^{i-1},f'_0(\tet,z,i)\dq (i-1)\tet^{i-2}z^{i-1},\\
 N_1(w,z,i)&=(-z^2+z^2wi+2z-zwi-1+zw)zwz^{i-1},D_1(w,z,i)\dq \frac{N_1(w,z,i)}{h(w,z)}\\
 N_2(w,z,i)&= -z^{i-1}zw(iz^2-2iz+2z+i-1-z^2+z^2wi-zwi+zw),D_2(w,z,i)\dq \frac{N_2(w,z,i)}{h(w,z)},\\
 C_2(j)&\dq -\frac{N_1(1,\ro,j)}{ h_w(1,\ro)}=j(-2+\srf)\ro^{j-1},\sum_j C_2(j)=C_2.
\end{align*}
 \[\mbox{Note that }C_2(j) \mbox{ is exponentially decreasing with }j. \mbox{ This will be the case for all following polyominoes. }\]
\[\FI(w,\tet,z,i)\dq w\tet z\lb f_0(\tet,z,i)+f_1(\tet,z) \frac{N_1(w,z,i)}{h(w,z)}+f_2(\tet,z)\frac{ N_2(w,z,i)}{h(w,z)}\rb.\]
 
But we also check that $G(\tet)$ is independent of $i$ (which is probabilistically obvious): equ.(\ref{E4}) is still valid and shows that $\frac{N_2(1,\ro,i)}{N_1(1,\ro,i)}$ is independent of $i$.
\subsection{The Markov chain (MC)}
We consider two successive columns $m_1,m_1+1$, of size $k,j$, such that their distances from the first and the last column are of order $\BO(n)$. Let $\Xi(m,n,m_1,k,j)$ denote the total number of dcc   with area $n$, width $m$,  column $m_1$ of size $k$,
 column $m_1+1$ of size $j$ and set $m_2:=m-m_1$. Theorem \ref{T22} leads to 
\bals
\Xi(m,n,m_1,k,j)&\dq \sum_{n_1} T(m_1,n_1,k,.)U(k,j)T(m_2,n-n_1,.,j),\\
 \frac{\Xi(m,n,m_1,k,j)}{T(m,n)}&\ds \sum_{n_1} \lb\frac{C_2}{\ro^{n_1}}\frac{e^{-(n_1-m_1\mu_2)^2/(2m_1\sig_2^2)}}{\sqrt{2\pi m_1}\sig_2}\pi(k)U(k,j) \frac{C_2(j)}{\ro^{n-n_1}}\frac{e^{-(n-n_1-m_2\mu_2)^2/(2m_2\sig_2^2)}}{\sqrt{2\pi m_2}\sig_2}\rb\\
 &\left/\lb\frac{C_2}{\ro^n}\frac{e^{-(n-m\mu_2)^2/(2m\sig_2^2)}}{\sqrt{2\pi m}\sig_2}\rb\right..
\end{align*}
Hence the asymptotic stationary distribution of two intermediate successive columns is given by
\[P(k,j)\dq \pi(k)U(k,j)C_2(j)= \pi(k)k C_2(j).\]
This leads to the asymptotic stationary  distribution  $\pi_2(k)$  and to the MC transition matrix $\Pi(k,j)$:
\bth
\bal
 \pi_2(k)&\dq \sum_j P(k,j)\dq \sum_j \pi(k)U(k,j)C_2(j),\mbox{  or }   \non\\
\pi_2(j)&\dq \sum_kP(k,j)\dq \sum_k \pi(k)U(k,j)C_2(j),  \non\\
\Pi(k,j)&\dq \frac{P(k,j)}{\sum_j P(k,j)}\dq \frac{U(k,j)C_2(j)}{\sum_j U(k,j)C_2(j)}= \frac{C_2(j)}{\sum_j C_2(j)}=\frac{C_2(j)}{C_2}
= j\ro^j,  \non\\
&\sum_k \pi_2(k) \Pi(k,j)\dq  \sum_k \sum_u   P(k,u) \frac{P(k,j)}{\sum_j P(k,j)} \dq  \pi_2(j),  \label{E41}\\
\pi_2(k)&= \pi(k)kC_2= k\ro^k.  \non
\end{align}
\ethGL
Eq. (\ref{E41})  confirms that $\pi_2(k)$  is the stationary distribution of $\Pi(k,j)$. 
This shows that the thickness of the polyomino is $\BO(1)$. This will be the case for all following polyominoes. 

A question we could ask is: is the chain reversible, i.e. is the following relation true?
\[\pi_2(k)\Pi(k,j)\dq \pi_2(j)\Pi(j,k)\equiv P(k,j)\dq P(j,k).\]
This  is satisfied here.

 The chain is irreducible and ergodic. Moreover, it is clear that the successive columns are independent and identically distributed (iid). Mean and variance of the stationary distribution $\pi_2(k)$ are given by
\[\mu_2\dq \sum_j \pi_2(j)j= \srf,\sig_x^2\dq \sum_j \pi_2(j)j^2-\mu_2^2= 2 \equiv \sig_2^2,
\mbox{ by independence of successive columns}.\]
We recover Sec.\ref{S21} results.
Several other  interesting relations can be derived. We have

\bals
 &\frac{C_2(k)}{\ro^{n}}\frac{e^{-(n-m\mu_2)^2/(2m\sig_2^2)}}{\sqrt{2\pi m}\sig_2}\ds
 \sum_j U(k,j)\frac{C_2(j)}{\ro^{n-k}}\frac{e^{-(n-k-(m-1)\mu_2)^2/(2(m-1)\sig_2^2)}}{\sqrt{2\pi (m-1)}\sig_2},\mbox{ hence }\\
C_2(k)&\dq \sum_j U(k,j)C_2(j)\ro^k,\\
\frac{\pi_2(k)}{\pi(k)}&\dq \sum_j U(k,j)C_2(j)\dq \frac{C_2(k)}{\ro^k},\\
\Pi(k,j)&\dq \frac{\pi(k)U(k,j)C_2(j)}{\pi_2(k)}\dq \frac{\ro^k}{C_2(k)}U(k,j)C_2(j)
\dq \frac{\pi(k)}{\pi_2(k)}U(k,j)\ro^j\frac{\pi_2(j)}{\pi(j)}.
\end{align*}
$\Pi(k,j)$ decreases exponentially with $j$ as well as $\pi_2(k)$ and $[\Pi]^\ell(k,j) $  converges exponentially fast to $\pi_2(j)$. The process is  $\Fi$-mixing (see 
Billingsley \cite[p.168ff]{BI68} and Appendix \ref{A1}). 
We will have the same properties for the following polyominoes. Also
\bals
\frac{\pi(j)}{\ro^j}&\dq \frac{\pi_2(j)}{C_2(j)}\dq \sum_k \pi(k)U(k,j),\\
C_2&\dq \sum  C_2(k) \mbox{ if starting with any  first column's size}.
 \end{align*}

\subsection{The perimeter conditioned on the width  $m$}                                                        \label{S13}
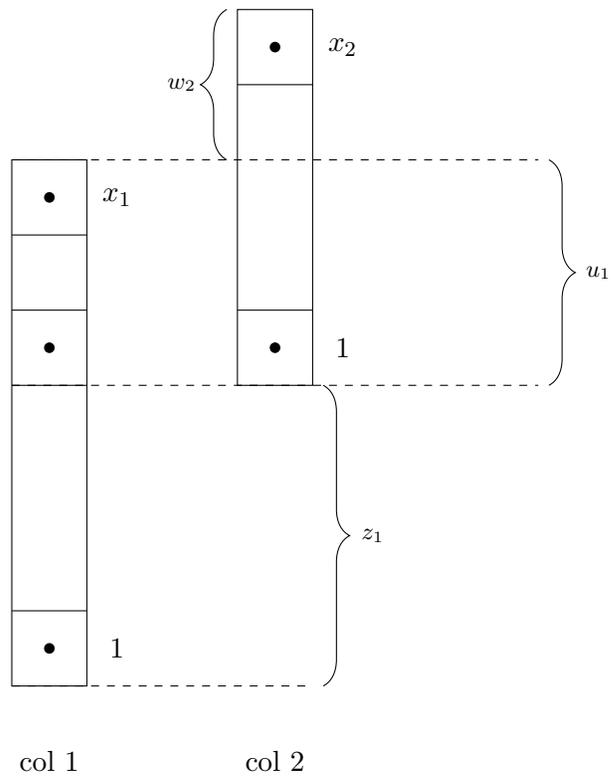
\begin{figure}[htbp]
\begin{center}
\begin{tikzpicture}[set style={{help lines}+=[dashed]}]
 \draw (2,0) rectangle (1,5);
 \draw (1,4) --  (2,4);
 \node[circle,fill=black!100,inner sep=0.05cm](s1)at(1.5,4.5){};
 \node at(2.4,4.5){$x_2$};
 \draw (1,1) --  (2,1);
 \node[circle,fill=black!100,inner sep=0.05cm](s1)at(1.5,0.5){};
 \node at(2.4,0.5){$1$};
 \draw (-2,-4) rectangle (-1,3);
 \draw (-2,2) --  (-1,2);
 \node[circle,fill=black!100,inner sep=0.05cm](s1)at(-1.5,2.5){};
 \node at(-0.6,2.5){$x_1$};
 \draw (-2,1) --  (-1,1);
 \draw (-2,0) --  (-1,0);
 \node[circle,fill=black!100,inner sep=0.05cm](s1)at(-1.5,0.5){};
 \draw (-2,-3) --  (-1,-3);
 \node[circle,fill=black!100,inner sep=0.05cm](s1)at(-1.5,-3.5){};
 \node at(-0.6,-3.5){$1$};
  \draw[dashed] (-2,3) --  (5,3);
  \draw[dashed] (-2,0) --  (5,0);
  \draw[dashed] (-2,-4) --  (2,-4);
\draw [decorate,decoration={brace,amplitude=10pt},xshift=-4pt,yshift=0pt]
(1,3) -- (1,5.0) node [black,midway,xshift=-0.6cm] 
{\footnotesize $w_2$};
\draw [decorate,decoration={brace,amplitude=10pt,mirror,raise=4pt},yshift=0pt]
(5,0.0) -- (5,3) node [black,midway,xshift=0.8cm] {\footnotesize
$u_1$};
\draw [decorate,decoration={brace,amplitude=10pt,mirror,raise=4pt},yshift=0pt]
(2,-4) -- (2,0) node [black,midway,xshift=0.8cm] {\footnotesize
$z_1$};
 \node at(-1.5,-5){col 1};
 \node at(1.5,-5){col 2};
\end{tikzpicture}
\end{center}
\caption{Two columns of a dcc polyomino and their related parameters.}
\label{F5}
\end{figure}

In the sequel, polyominoes of width $m$  can be seen as a sequence of id RV (columns) $x_i,i=1,2,\ldots,x_m$  by the MC   $\Pi(k,j)$.
In this section, we fix  $m$ and analyze the asymptotic properties of the perimeter $P_m$. For dcc, we have the following 
notations and probabilistic relations: see Fig.\ref{F5},

\bals
 w_d&:=(\mbox{ upper cell position of column }d)-(\mbox{ upper cell position of column }d-1),\\
 z_d&:=(\mbox{ lower cell position of column }d+1)-(\mbox{ lower cell position of column }d),\\
z_d \geq 0 & \mbox{ depends only on }x_d \mbox{ and is uniformly distributed } (0,x_{d}-1),\\
u_d \geq 1 & \mbox{ depends only on }x_d \mbox{ and is uniformly distributed } (1,x_d),\\
w_d&=x_d-u_{d-1},z_d=x_d-u_d,T_d:=|w_d|+z_d.
\end{align*}
 the total perimeter  $P_m$  the asymptotic  vertical perimeter $Q_m$
 the asymptotic  total  perimeter  $R_m$
are  given by  
\bals
P_m &\dq Q_m+x_0+x_m+2m,Q_m\dq \sum_1^m T_d,R_m\dq Q_m+2m,X_m\dq \sum_1^m x_d,\\
\V(X_m)&\ds m\sig_2^2,(\V(.) \mbox{ denoting the variance)},\\
\E(P_m)&\ds \E(Q_m)+2m,\V(P_m)\ds \V(Q_m),\E(X_m)\ds  m \mu_2,Cov(P_m,X_m)\ds Cov(Q_m,X_m).
\end{align*}
We compute now the first probability densities
\bals
P_u(i)&:=\P(u_1=i)=\sum_{k=i}^\II \P(x_1=k)\frac1k=\sum_{k=i}^\II\frac{\pi_2(k)}{k}= \ro^i\frac{\srf+1}{2},\\
f(t)&:=\P(|w|)=t)=\sum_{r=1}^\II P_u(r)\pi_2(r+t)+\sum_{r=t+1}^\II P_u(r)\pi_2(r-t),t>0,\\
f(0)&=\sum_{r=1}^\II P_u(r)\pi_2(r)= \frac{\srf+1}{10}.
\end{align*}

The successive moments we need are computed as follows. The mean and variance of $T_d$ will be denoted by $\mu_3, \sig^2_3$,

\bals
E_w&:=\E(|w|)=\sum_1^\II f(t)t= \frac{7\srf-3}{10},\\
E_z&=\sum_{j=1}^\II \pi_2(j)\sum_{s=0}^{j-1}\frac{s}{j}= \frac{\srf-1}{2},\\
\mu_3&\dq \E(T_d)=E_w+E_z= \frac{6\srf -4}{5},\E(Q_m)\ds m\mu_3.\\
\mbox{Let}\\
S_1&=\E(w_2^2)=\sum_1^\II f(t)t^2,\\
S_2&=\E(|w_2|z_2)=\sum_{r=1}^\II P_u(r)\lb \sum_{j=r}^\II \pi_2(j)(j-r)\frac1j \sum_{s=0}^{j-1}s 
+\sum_{j=1}^{r-1} \pi_2(j)(r-j)\frac1j \sum_{s=0}^{j-1}s\rb ,\\
S_3&=\E(z_2^2)=\sum_{j=1}^\II \pi_2(j)\frac1j \sum_{s=0}^{j-1}s^2 ,\mbox{ then }\\
\E(T^2_2)&=\E([|w_2|+z_2]^2)=\E([|x_2-u_1|+z_2]^2)= S_1+2S_2+S_3=-\frac{52\srf}{25}+\frac{62}{5},\\
\sig^2_3&\dq \E(T^2_2)-\mu_3^2=\frac{-4\srf+114}{25}.
\end{align*}
Let
\bals
S_4&=\E(|w_2||w_3|)=\E(|x_2-u_1||x_3-u_2|)\\
&=\sum_{r=1}^\II P_u(r)
\lb \sum_{j=r}^\II \pi_2(j)(j-r)\frac1j\sum_{v=1}^j\lb \sum_{\ell=v}^\II\pi_2(\ell)(\ell-v) +\sum_{\ell=1}^{v-1}\pi_2(\ell)(v-\ell)\rb\right.\\
 &+\left. \sum_{j=1}^{r-1} \pi_2(j)(r-j)\frac1j\sum_{v=1}^j\lb \sum_{\ell=v}^\II\pi_2(\ell)(\ell-v) +\sum_{\ell=1}^{v-1}\pi_2(\ell)(v-\ell)\rb\rb,\\
S_5&=\E(|w_2|z_3)=\E(|x_2-u_1|z_3)=E_w E_z,\\
S_6&=\E(z_2z_3)=E_z^2,\\
S_7&=\E(z_2|w_3|)=\sum_{j=1}^\II \pi_2(j)\frac1j\sum_{v=1}^j (j-v)\lb\sum_{\ell=v}^\II\pi_2(\ell)(\ell-v) +\sum_{\ell=1}^{v-1}\pi_2(\ell)(v-\ell)\rb,\mbox{ then }\\
E(T_2 T_3)&=S_4+S_5+S_6+S_7= \frac{157}{20}-\frac{91\srf}{50}.
\end{align*}
Finally, we obtain
\bals
\V(Q_m)&\ds m \sig^2_Q, \mbox{ with}\\
\sig^2_Q&=\sig^2_3+2(E(T_2T_3)-\mu_3^2)= \frac{229}{50}+\frac{\srf}{25}.
\end{align*}
The covariance and correlation coefficient are computed as follows. Let

\bals
S_8&=\E[x_2(|w_2|+z_2)]=\\
&\sum_{r=1}^\II P_u(r)\lb\sum_{j=r}^\II \pi_2(j)(j-r)j+\sum_{j=1}^{r-1} \pi_2(j)(r-j)j\rb
+\sum_{j=1}^\II \pi_2(j)\frac1j \sum_{s=0}^{j-1}js,\\
S_9&=\E[x_2(|w_3|+z_3)]\\
&=\sum_{j=1}^\II \pi_2(j)j\frac1j\sum_{v=1}^j \lb\sum_{\ell=v}^\II\pi_2(\ell)(\ell-v) +\sum_{\ell=1}^{v-1}\pi_2(\ell)(v-\ell)\rb+\mu_2 E_z,\mbox{ then }\\
Cov(X_mQ_m)&\ds mC_{X,Q},C_{X,Q}=S_8+S_9-2\mu_2\mu_3= \frac{13}{5}+\frac{\srf}{25},\\
\ro(X_m,X_m)&\ds \ro_{X,Q},\ro_{X,Q}\dq \frac{C_{X,Q}}{\sig_2\sig_Q}= \frac{ (65+\srf)\sqrt{2}}{5(458+4\srf)^{1/2}}.
 \end{align*}
\subsection{The perimeter conditioned on the area $n$}
In this subsection, we obtain  convergence to Brownian motions of $X(\lf \nu t \rf),Q(\lf \nu t \rf)$, and we compute asymptotic mean and variance of $R$ conditioned on $n$, denoted by $\Rb(n)$. First of all we fix the width to $m$. 
By the function central limit theorem for dependent random variables (RV) and by the mixing property (see, for instance, Billingsley \cite[p.168ff]{BI68}),  we obtain
the following conditioned on $m$ convergences (where $B_.(t)$ are standard Brownian motions), and these convergences will be valid for the following polyminoes,
\bth
\bals
&\frac{X(\lf \nu t \rf)-\mu_2\nu t}{\sqrt{\nu}\sig_2}\conv  B_1(t),\nu\ra\II,t\in[0,1],\\
&\frac{Q(\lf \nu t \rf)-\mu_3\nu t}{\sqrt{\nu}\sig_Q}\conv B_2(t),\nu\ra\II,t\in[0,1],\\
\end{align*}
\ethGL

 The convergence of $Q$ is due to the fact that $T_d$  depends only on $x_d,x_{d-1}$
 and the mixing property still holds.  Therefore 
\bals
X(m)&\ds  m\mu_2+\sqrt{m}\sig_2\ta_1,\ta_1\dq \BN(0,1),\\
Q(m)&\ds  m\mu_3+\sqrt{m}\sig_Q\ta_2,\ta_2\dq \BN(0,1),\ta_2\dq \ro_{X,Q}\ta_1+\sqrt{1-\ro_{X,Q}^2}\ta_3,\ta_3\dq \BN(0,1),\\
R(m)&\dq Q(m)+2m,\mbox{ hence }\\
&\frac{R(m)-m(\mu_3+2)}{\sqrt{m}\sig_Q}\ds \ro_{X,Q}\frac{n-m\mu_2}{\sqrt{m}\sig_2}+\sqrt{1-\ro_{X,Q}^2}\ta_3,\\
\E(R(m))&\ds n\al+m\be,\al\dq \frac{\sig_Q\ro_{X,Q}}{\sig_2},\be\dq \mu_3+2-\ro_{X,Q}\mu_2\frac{\sig_Q}{\sig_2},\\
 \V(R(m))&\ds m\gam,\gam\dq \sig_Q^2(1-\ro_{X,Q}^2).
 \end{align*}
By equ.(\ref{E3}), we derive
\bals
\E\lp e^{ \ii\tet \Rb(n)}\rp&\ds \E\lb e^{\ii\tet[\al n+\be m]-\frac{\tet^2}{2}\gam m}  \rb\\
&\ds \E\lb e^{\ii\tet[\al n+\be(n\mu_1+\sqrt{n}\sig_1\ta_0)]-\frac{\tet^2}{2}\gam (n\mu_1+\sqrt{n}\sig_1\ta_0)}  \rb\\
&\ds e^{\ii\tet[\al n+\be n\mu_1]-\frac{\tet^2}{2}\gam n\mu_1
-\frac{1}{2}\lb\tet \be \sqrt{n}\sig_1+\ii \frac{\tet^2}{2}\gam \sqrt{n}\sig_1\rb^2},
\end{align*}
and setting $\tet\dq \tet/\sqrt{n}$, we finally obtain the following result
\bth
\bals
\frac{\Rb(\lf n t \rf)-\mu_4 n t}{\sqrt{n}\sig_4}&\conv B_4(t),n\ra\II,t\in[0,1],\\
 \frac{\Rb(n)-n\mu_4}{\sqrt{n}\sig_4}&\ds \ta_4,\ta_4\dq \BN(0,1),\\
\mu_4&\dq \al+\be\mu_1\dq \mu_1(\mu_3+2)=\frac65+\frac{6\sqrt{5}}{25}=1.736656315\ldots,\\
\sig_4^2&\dq \gam\mu_1+\sig_1^2\be^2=\frac{17\srf}{50}-\frac{19}{125}=.6082631123\ldots
\end{align*}
\ethGL
Note that, in some cases (here and in the wall polyomino case) our technique leads to exact values for $\mu_4,\sig_4^2$. 
We have made extensive simulations to check our results. We first  construct $T_s=400$ times a MC $X(d),d=1..m$ based on  $\Pi(k,j)$ with $m=400$ steps. This allows to check the values $\mu_2,\sig_2^2$. The fit is excellent. Next we extend ( or contract) each MC  such that  $X(m^*)=n$,    with $n= \lf m/\mu_1\rf$.  Based on each MC, we build  on each column $d$ a vertical perimeter $T_d$ following the distribution of Sec. \ref{S13}. We then compare the observed distribution of the vertical perimeter $Q(m^*)$ with the theoretical parameters $\mu_4^*,\sig_4^{*2}$, where
\[\mu_4^*\dq  \mu_1\mu_3,\sig_4^{*2}\dq \gam\mu_1+\sig_1^2\be^{*2},\be^* \dq \mu_3-\ro_{X,Q}\mu_2\frac{\sig_Q}{\sig_2}.\]
Indeed, we don't take the $2m$ horizontal perimeter into account. The fit is quite good.

Let us finally make four  remarks
\bit
\item Actually, we can use Bender's theorems in another way: it is possible to derive large deviations results for all our convergence to Gaussian variables theorems: see for instance Louchard \cite{GL297}. We will not detail these applications here.
\item If we compare exact distributions with the Gaussian limits, we observe a bias, for instance for $n=30$. This can be corrected with Hwang  \cite [Thm.2]{HW94}, Hwang \cite{HW99}. (See also Flajolet and Sedgewick \cite [Lemma IX.1]{FLSe09}). See an example in Louchard \cite{GL297}.
\item The maximum thickness of the polyomino can also be analyzed. See for instance Louchard \cite{GL299}.
\item The trajectories of the polyomino (upper and lower trajectories) lead themselves to Brownian motions: see \cite{GL96} for dcc and \cite{GL299} for cc, dc. 
The polyomino can be seen as a Brownian motion with some thickness. A more detailed analysis will be the object of a future report. 
\eit
\subsection{A comparison with known GF}
In some cases, we know
 the joint GF of $n$ and $R_n$. For instance, for dcc, this is given  in Bousquet-M\'elou \cite[(10)]{BM96}
\footnote{ see Appendix \ref{A2} }. If we set $q=z,y=1,x=v$ 
in the denominator of (10), we  we recover of course $h(w,z)$, leading to the root $\ro$. If we set $q=z,y=x=v$,    we obtain
\footnote{We use the Pochhammer symbol: $(a;z)_j:=(1-a)(1-az)\ldots(1-az^{j-1})$}
\[F(v,z)=L_0(1)=1-\sum_{j=1}^\II \frac{v^j(v-1)^{j-1}z^{j(j+1)/2}}{(z;z)_j(vz;z)_{j-1}(vz;z)_j}.\]
This corresponds to the half-perimeter $R_n/2$.   Using again Bender's theorems,  we obtain
\[\frac{R_n/2-\mub_4 n}{\sigb_4\srn } \ds \BN(0,1),\mbox{ and }  2\mub_4=1.736656315\ldots,4\sigb_4^2=0.6082631120\ldots,\]
which fits with $\mu_4,\sig_4^2$.
\section{The cc perimeter}
\subsection{The generating functions}
A column-convex polyomino (cc), is made of contiguous columns such that at least one cell of each column must be adjacent to some cell of the previous column. 
We first recall from \cite{GL297} the main expressions we need: starting with any first column's  size, and with the same notations as in the previous section, 
\bals
 U(k,j)&=k+j-1,\tilde{f}_1(\tet,z)= \tet^2 z^2/(1-\tet z)^2,\tilde{f}_2(\tet,z)= \tet z/(1-\tet z),
\tilde{f}_0(\tet,z)= \tilde{f}_2(\tet,z),\\ 
f_1(\tet,z)&= 1/(1-\tet z)^2,f_2(\tet,z)= 1/(1-\tet z),f_0(\tet,z)= f_2(\tet,z),\\
 N_1(w,z)&= wz(z-1)^3,N_2(w,z)= -z^2w(z^2-2z+1+zw),\\
 h(w,z)&= z^4(w-1)+z^3(w^2-w+4)-z^2(w+6)+z(w+4)-1, \\
&\mbox{ see also  Flajolet and Sedgewick \cite[p.366]{FLSe09} }.\\
&\mbox{ Solving }4\ro^3-7\ro^2+5\ro-1=0,\mbox{ this gives }\\
 \ro&= -C_3^{1/3}+\frac{11}{144C_3^{1/3}}+\frac{7}{12},C_3= \frac{71}{1728}+\frac{\sqrt{177}}{288},\\
 G'(1)&= \frac{1-2\ro}{\ro(1-\ro)},\pi(j)= \ro^j\lp \frac{1-3\ro+\ro^2}{\ro(1-\ro)}+j\rp,\\
 \mu_1&= \frac{11\ro^2-9\ro+5}{4(12\ro2-14\ro+5)},\mu_2=1/\mu_1,\sig_1^2=\frac{-1478891+6578899\ro-5346249\ro^2}{256(-26895+104919\ro-44437\ro^2)},\\
 C_2&= \frac{5-13\ro+7\ro^2}{11-35\ro+41\ro^2},\\
\sig_2^2&\dq \frac{\sig_1^2}{\mu_1^3}=-16\frac{4579\ro^2-9681\ro+1753}{28283\ro^2-57561\ro+24097}.\\
&\mbox{Starting with a column of size }i,\\
 N_1(w,z,i)&=  wz^i(z-1)^2(z^2wi-z^2+zw-zwi+2z-1),\\ 
N_2(w,z,i)&=\\
& (z-1)(-iz^3+z^3-3z^2+z^2w+3iz^2+z^2wi-zwi-3iz+3z+zw+i-1)zwz^i,\\
 C_2(j)&=(aj+b)\ro^j,a:=\frac{5-13\ro+7\ro^2}{11-35\ro+41\ro^2} ,b:=\frac{3-11\ro+17\ro^2}{11-35\ro+41\ro^2},\sum_1^\II C_2(j)=C_2.\\
&\mbox{We have two relations:}\\
 \frac{b}{a}&=(1-3\ro+\ro^2)/(\ro(1-\ro)),\\
&\sum_{j=1}^\II (k+j-1)\ro^j(aj+b)=ak+b.
\end{align*} , 
\subsection{The Markov chain}
We have here
\bals
\pi(j)&= \ro^j(aj+b)/a=C_2(j)/a,U(k,j)=k+j-1,\\
P(k,j)&\dq \pi(k)U(k,j)C_2(j)=\ro^k(ak+b)(k+j-1)(aj+b)\ro^j/a.\\
P(k,j)&\mbox{ is symmetric, hence the chain is reversible },\\
\Pi(k,j)&= \frac{1}{ak+b}(k+j-1)(aj+b)\ro^j\\
&=\frac{\ro^k}{\pi(k)}(k+j-1)\pi(j)=\frac{\ro^k}{C_2(k)}(k+j-1)C_2(j),
\pi_2(k)\dq  (ak+b)^2\ro^k/a,\\
\mu_2&\dq \sum \pi_2(k) k=  \frac{ 16(4131-14923z+14001z^2)}{5185-7673z+33755z^2}\equiv \frac{1}{\mu_1},\\
\sig_x^2&\dq \sum \pi_2(k) k^2-\mu_2^2\\
& = \frac{50790154312925840-247592802999061008z+268780110914590000z^2}
{-5468815736218009+23664229539220113z-21285767156057123z^2}.\\
 \end{align*}
\subsection{The perimeter conditioned on $m$ and $n$}
\begin{figure}[htbp]
\begin{center}
\begin{tikzpicture}[set style={{help lines}+=[dashed]}]
 \draw (2,-7) rectangle (1,5);
 \draw (1,4) --  (2,4);
 \node[circle,fill=black!100,inner sep=0.05cm](s1)at(1.5,4.5){};
 \draw (1,1) --  (2,1);
 \node[circle,fill=black!100,inner sep=0.05cm](s1)at(1.5,0.5){};
 \draw (1,2) --  (2,2);
 \draw (1,3) --  (2,3);
 \node[circle,fill=black!100,inner sep=0.05cm](s1)at(1.5,2.5){};
 \draw (-2,-4) rectangle (-1,3);
 \draw (-2,2) --  (-1,2);
 \node[circle,fill=black!100,inner sep=0.05cm](s1)at(-1.5,2.5){};
 \draw (1,-6) --  (2,-6);
 \node[circle,fill=black!100,inner sep=0.05cm](s1)at(1.5,-6.5){};
 \draw (-2,1) --  (-1,1);
 \draw (-2,0) --  (-1,0);
 \node[circle,fill=black!100,inner sep=0.05cm](s1)at(-1.5,0.5){};
 \draw (-2,-3) --  (-1,-3);
 \draw (1,-3) --  (2,-3);
 \node[circle,fill=black!100,inner sep=0.05cm](s1)at(-1.5,-3.5){};
 \node[circle,fill=black!100,inner sep=0.05cm](s1)at(1.5,-3.5){};
  \draw[dashed] (-1,3) --  (2,3);
  \draw[dashed] (-1,2) --  (2,2);
  \draw[dashed] (-1,0) --  (2,0);
  \draw[dashed] (-1,1) --  (2,1);
  \draw[dashed] (-1,-4) --  (2,-4);
  \draw[dashed] (-2,-7) --  (2,-7);
  \draw[dashed] (-2,5) --  (2,5);
\node at(0.8,-6.5){};
 \draw [decorate,decoration={brace,amplitude=10pt},xshift=-4pt,yshift=0pt]
(-2,3) -- (-2,5.0) node [black,midway,xshift=-0.6cm] 
{\footnotesize $w_v$};
\draw [decorate,decoration={brace,amplitude=10pt},xshift=-4pt,yshift=0pt]
(-2,-4) -- (-2,3) node [black,midway,xshift=-0.6cm] 
{\footnotesize $k$};
\draw [decorate,decoration={brace,amplitude=10pt},xshift=-4pt,yshift=0pt]
(-2,-7) -- (-2,-4) node [black,midway,xshift=-0.6cm] 
{\footnotesize $z_v$};
\draw [decorate,decoration={brace,amplitude=10pt,mirror,raise=4pt},yshift=0pt]
(2,-7) -- (2,5) node [black,midway,xshift=0.8cm] {\footnotesize
$j$};
 \node at(-1.5,-8){col $v$-1};
 \node at(1.5,-8){col $v$};
\end{tikzpicture}
\end{center}
\caption{Two columns of a cc polyomino and their related parameters.}
\label{F6}
\end{figure}
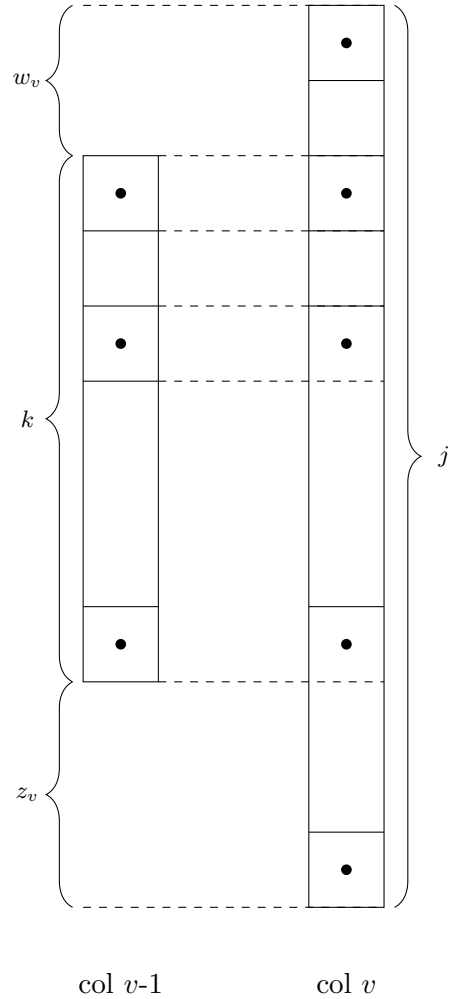

 For cc, we have the following 
notations and probabilistic relations: see Fig.\ref{F6},  

\bals
 w_v&:=(\mbox{ upper cell position of column }v-1)-(\mbox{ upper cell position of column }v),\\
 z_v&:=(\mbox{ lower cell position of column }v-1)-(\mbox{ lower cell position of column }v),\\
w_v,z_v,  & \mbox{ depend only on }k=x_{v-1} \mbox{ and } j=x_{v} ,\\
z_v &\mbox{  is uniformly distributed } (-(k-1),j-1),\\
w_v &\mbox{  is uniformly distributed } (-(j-1),k-1),w_v=k-j+z_v,T_v= |w_v|+|z_v|.\\
&\mbox{ The moments are computed as follows }\\
E_w(k,j)&=E_z(k,j)=\E(\left.\rule{0mm}{4mm} |z|\right|k,j)\dq \frac{1}{k+j-1}\lb  \sum_{i=-(k-1)}^0 (-i)+\sum_{i=1}^{j-1} (+i) \rb\dq
\frac{k^2+j^2-k-j}{2(j+k-1)},\\
E_{w^2}(k,j)&=E_{z^2}(k,j)=\E(\left.\rule{0mm}{4mm} z^2\right|k,j)\dq \frac{1}{k+j-1}\lb  \sum_{i=-(k-1)}^{j-1}  i^2 \rb
\dq 1/3j^2-1/6j-1/3jk-1/6k+1/3k^2,\\
E_{zw}(k,j)&=\E(\left. \rule{0mm}{4mm} |z||w|\right|k,j).\\
&\mbox{ For the case }j>k,\mbox{ we set }E^p,A= k-j+i\\
E_{zw}^p(k,j)&=\frac{1}{k+j-1}\lb   \sum_{i=-(k-1)}^0 (-i)(-A)+\sum_{i=j-k}^{j-1} iA+\sum_{i=1}^{j-k-1} i(-A )\rb\\
&= 1/6(-j+3k+j^3-6jk-3k^3-3kj^2+9jk^2)/(j+k-1).\\
&\mbox{ For the case }j\leq k,\mbox{ we set }E^m,A= k-j+i, \\
&m \mbox{ is here related to the case }j\leq k,\mbox{ and not to the width }\\
E_{zw}^m(k,j)&= \frac{1}{k+j-1}\lb   \sum_{i=j-k}^0 (-i)A+\sum_{i=-(k-1)}^{j-k-1} (-i)(-A)+\sum_{i=1}^{j-1} iA\rb\\
&= -1/6(-3j+k+3j^3+6jk-k^3-9kj^2+3jk^2)/(j+k-1),\\
&\mbox{note that we have some symmetry here: }E^m_{.}(k,j)\equiv E^p_{.}(j,k).\\
&\mbox{ For } y \mbox{ denoting   a random variable depending on } k,j \mbox{ with mean }E_y(k,j),\mbox{ we set }\\
\Eb(y)\dq & \sum_k\pi_2(k)\sum_j\Pi(k,j)E_y(k,j),\\
&\mbox{ for } f(k,j) \mbox{ denoting  a function  depending on } k,j ,
\mbox{ we set }\Eb(f)\dq  \sum_k\pi_2(k)\sum_j \Pi(k,j)f(k,j),\\
&\mbox{ we compute }\\
\mu_3&\dq \E(T_d)=2\Eb(z)=1.962459470\ldots,\\
\sig_3^2&\dq \sig^2(T_d)=\Eb\lp \lb |w|+|z|-2\Eb(z)\rb^2\rp =\Eb\lp\lb |w|+|z|-2E_z+2E_z-2\Eb(z)\rb ^2\rp=S_1+S_2,\\
S_1&= \Eb\lp 2E_{z^2}+2E_{zw}-(2E_z)^2  \rp,\\
S_2&= 4\lp \Eb(E_z^2-\Eb(z)^2\rp,\\
\sig_3^2&=\Eb\lp (2E_{z^2}+2E_{zw}\rp-\mu_3^2=2.387549945\ldots,\\
\Eb(zw)&=\sum_{k=1}^\II \pi_2(k)\sum_{j=k+1}^\II\Pi(k,j)E_{zw}^p(k,j)+\sum_{k=1}^\II\pi_2(k)\sum_{j=1}^k\Pi(k,j)E_{zw}^m(k,j),\\
\Eb(z^2)&=\sum_{k} \pi_2(k)\sum_{j}\Pi(k,j)E_{z^2}(k,j),\\
\end{align*}
\bals
\sig_Q^2&\dq \sig_3^2+2\sum_2^\II C_\ell,C_\ell\dq\E(T_1T_\ell)-\mu_3^2,T_k-\mu_3=[|z_k|+|w_k|-\mu_3].\\
&\mbox{ Generalizing some results from \cite{GL297}, and setting}\\
y_{1,i}&\dq \mbox{ random variable depending only on }x_{i-1}=k,x_{i}=j,\mbox{ with mean }\E(y_{1,i})\dq F_1(k,j),\\
y_{2,i}&\dq \mbox{ random variable depending only on }x_{i-1}=k,x_{i}=j,\mbox{ with mean }\E(y_{2,i})\dq F_2(k,j),\\
\mub_r&\dq \E(y_r)\dq \sum_u \pi_2(u)\sum_j \Pi(u,j) F_r(u,j),\\
&\mbox{ we derive}\\
\E(y_{1,1}.y_{2,2})& \dq \sum_u \pi_2(u)\sum_j \Pi(u,j) F_1(u,j) \sum_k  \Pi(j,k)F_2(j,k)\\
&\dq \sum_u\pi_2(u)\sum_j\Pi(u,j)F_1(u,j)\frac{1}{C_2(j)}\ro^j\sum_k U(j,k)C_2(k)F_2(j,k),\\
\E(y_{1,1}.y_{2,3})& \dq \sum_u\pi_2(u)\sum_j\Pi(u,j)F_1(u,j)\frac{1}{C_2(j)}\sum_\ell[\ro^{j+\ell}U(j,\ell)]
\sum_k U(\ell,k)C_2(k)F_2(\ell,k).\\
&\mbox{ The Markov property leads to compute }\\
&\sum_{\ell=2}^\II[\E(y_{1,1}.y_{2,\ell})-\mub_1\mub_2]\dq \\
\lim_{w\ra 1}\lb  \rule{0mm}{8mm}\sum_u\sum_j\sum_\ell\sum_k \right. & \left. \pi_2(u)\Pi(u,j)F_1(u,j)\frac{1}{C_2(j)}
\lb[\tet^\ell]\FI(w,\tet,\ro,j)\rb
  U(\ell,k)C_2(k)F_2(\ell,k)-\frac{w}{1-w} \mub_1\mub_2 \rule{0mm}{8mm} \rb.
	\end{align*}
	 Starting with a first column of size $ i$, we set \footnote{We use the Iverson Bracket, as advocated by D.E.Knuth: $\kl P \kr=1$ if $P =true,=0$ otherwize }
	\bals
M_1(j,\ell)&\dq [\tet^\ell]\tet^j\ro^j=\ro^j\kl j=\ell \kr,\\
\Xi_3(F_1,F_2)&\dq \sum_u\sum_j\sum_\ell\sum_k \pi_2(u)\Pi(u,j)F_1(u,j)\frac{1}{C_2(j)}M_1(j,\ell) U(\ell,k)C_2(k)F_2(\ell,k)\\
&\dq \sum_u\sum_j\sum_k \pi_2(u)\Pi(u,j)F_1(u,j)\frac{1}{C_2(j)}\ro^j U(j,k)C_2(k)F_2(j,k),\\
\Xi_4(F_1,F_2)&\dq \\
\lim_{w\ra 1}\lb \rule{0mm}{8mm}\sum_u\sum_j\sum_\ell\sum_k\right.& \left.  \pi_2(u)\Pi(u,j)F_1(u,j)\frac{1}{C_2(j)}
\lb[\tet^\ell]\frac{\Fi(w,\tet,\ro,j)}{h(w,\ro)}\rb  U(\ell,k)C_2(k)F_2(\ell,k)-\frac{1}{1-w} \mub_1\mub_2 \rule{0mm}{8mm}\rb\\
\Fi(w,\tet,\ro,j)&\dq \tet\ro[f_1(\tet,\ro) N_1(w,\ro,j)+f_2(\tet,\ro) N_2(w,\ro,j)].\\
&\mbox{Note that, in order to simplify the limits,  we have divided our expressions by }w. 
\end{align*}
\bals
&\mbox{Set now }w=1-\eps,
\mbox{this leads to }\\
\lim_{\eps\ra 0}\lb \rule{0mm}{8mm}\sum_u\sum_j\sum_\ell\sum_k\right.&   \pi_2(u)\Pi(u,j)F_1(u,j)\frac{1}{C_2(j)}\times\\
\times [\tet^\ell]&\left. \lb -\frac{\Fi(1,\tet,\ro,j)}{h_w(1,\ro)\eps}+\lb    \frac{\Fi_w(1,\tet,\ro,j)}{h_w(1,\ro)}
-\frac{\Fi(1,\tet,\ro,j)h_{ww}(1,\ro}{2h^2_w(1,\ro)} \rb \rb U(\ell,k)C_2(k)F_2(\ell,k)-\frac{\mub_1\mub_2}{\eps} \rule{0mm}{8mm}\rb.\\
\end{align*}
\bals
&\mbox{We must first dispense from the singularity. But }\\
\frac{\Fi(1,\tet,\ro,j)}{h_w(1,\ro)}&\dq \frac{N_1(1,\ro,j)G(\tet}{h_w(1,\ro)}\dq -C_2(j)G(\tet),[\tet^\ell](-C_2(j)G(\tet))
\dq -C_2(j)\pi(\ell),\mbox{ hence }\\
\sum_u\sum_j\sum_\ell\sum_k&   \pi_2(u)\Pi(u,j)F_1(u,j)\frac{1}{C_2(j)}\lb \frac{C_2(j)\pi(\ell)}{\eps}\rb              U(\ell,k)C_2(k)F_2(\ell,k)\\
&\dq \frac{1}{\eps}\sum_u\sum_j\sum_\ell\sum_k \pi_2(u)\Pi(u,j)F_1(u,j)\pi_2(\ell)\Pi(\ell,k)F_2(\ell,k)\dq \frac{\mub_1\mub_2}{\eps},\\
&\mbox{and the singularity is removed. Set}\\
M_2(j,\ell)&\dq [\tet^\ell]\lb    \frac{\Fi_w(1,\tet,\ro,j)}{h_w(1,\ro)}-\frac{\Fi(1,\tet,\ro,j)h_{ww}(1,\ro)}{2h^2_w(1,\ro)} \rb,\\
\Xi_4(F_1,F_2)&\dq \sum_u\sum_j\sum_\ell\sum_k \pi_2(u)\Pi(u,j)F_1(u,j)\frac{1}{C_2(j)} M_2(j,\ell) U(\ell,k)C_2(k)F_2(\ell,k).\\
\end{align*}

Finally
\[ \sum_{\ell=2}^\II[\E(y_{1,1}.y_{2,\ell})-\mub_1\mub_2]\dq \Xi_5(F_1,F_2)\dq \Xi_3(F_1,F_2)+\Xi_4(F_1,F_2).\]
This relation will also be used in some following polyominoes.\\
Note that, in our previous expressions, we can use
\[\pi_2(u)\Pi(u,j)/C_2(j)\dq \pi(j)U(u,j).\]
We now obtain
\bals
\sig_Q^2&=\sig_3^2+2\Xi_5(2E_z,2E_z)=3.8341042755\ldots,\\
\sig_X^2&=\sig_x^2+2\Xi_5(j,k)\equiv \sig_2^2,\mbox{ this has been explicitly checked by direct computation},\\
Cov(X(m),Q(m))&=\E[\sum_1^m(x_i-\mu_2)\cdot \sum_1^m[|z_k|+|w_k|-\mu_3]]\sim mC_{X,Q},C_{X,Q}=S_6+S_7,\\
S_6&=\sum_k\pi_2(k)\sum_j\Pi(k,j)j\cdot 2E_z(k,j)-\mu_2\mu_3,\\
S_7&=\Xi_5(j,2E_z)+\Xi_5(2E_z,k),\\
\ro_{X,Q}&\dq \frac{C_{X,Q}}{\sig_2\sig_Q}=0.8873927438\ldots,\mu_4=1.7952896266\ldots,\sig_4^2=0.4588988471\ldots\\
\end{align*}
Again, we have made extensive simulations to check our results. The fit is quite good.

\begin{figure}[htbp]
	\centering
		\includegraphics[width=0.8\textwidth,angle=0]{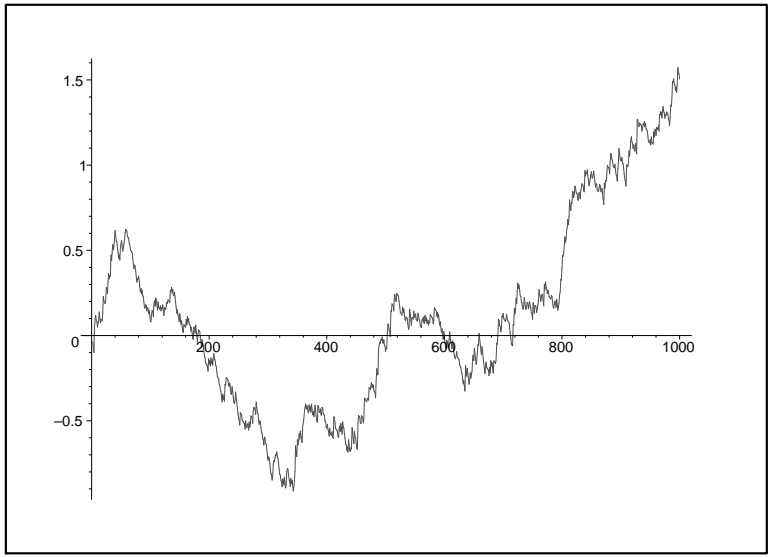}
	\caption{Simulation of 
$\frac{X(v)-\mu_2 v}{\sqrt{1000}\sig_2},v=1..1000$}
	\label{F1}
\end{figure}

\begin{figure}[htbp]
	\centering
		\includegraphics[width=0.8\textwidth,angle=0]{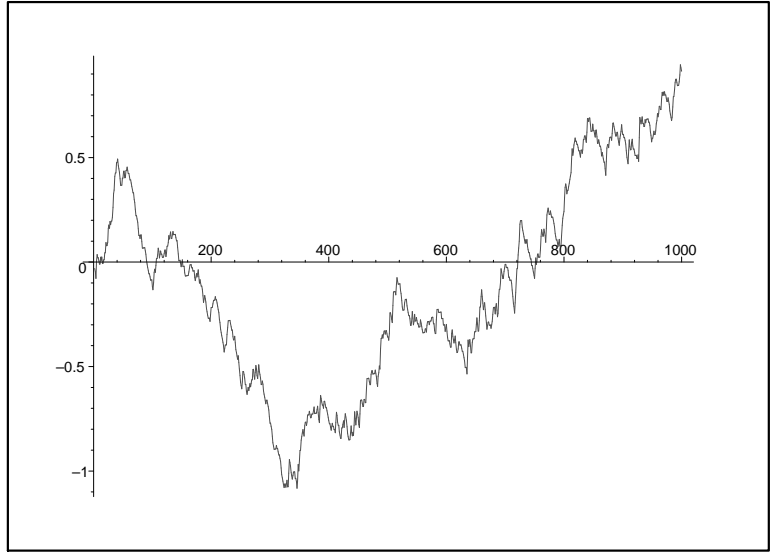}
	\caption{Simulation of 
$\frac{Q(v)-\mu_3 v}{\sqrt{1000}\sig_Q},v=1..1000$}
	\label{F2}
\end{figure}

\begin{figure}[htbp]
	\centering
		\includegraphics[width=0.8\textwidth,angle=0]{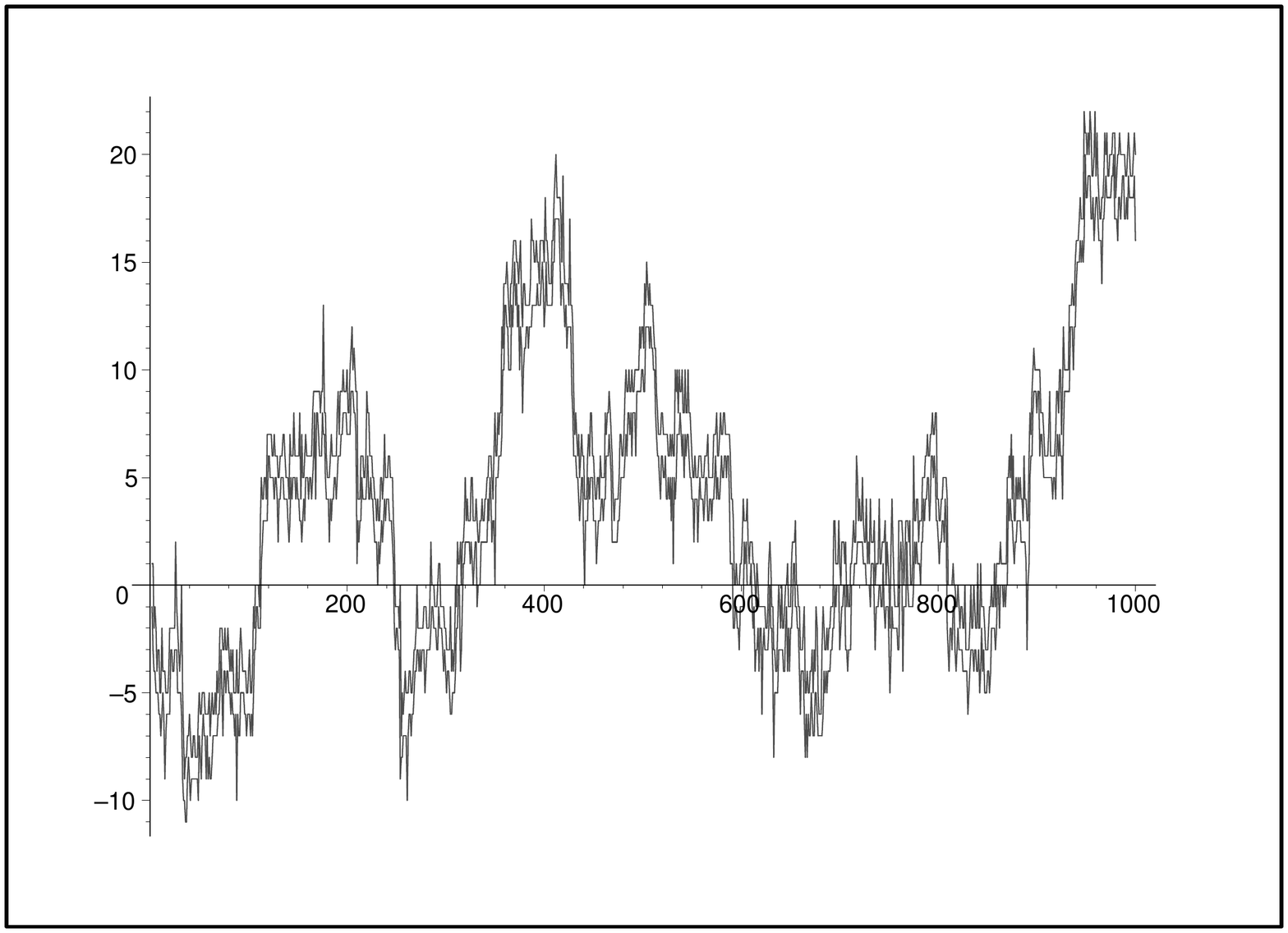}
	\caption{A typical polyomino}
	\label{F3}
\end{figure}

\begin{figure}[htbp]
	\centering
		\includegraphics[width=0.8\textwidth,angle=0]{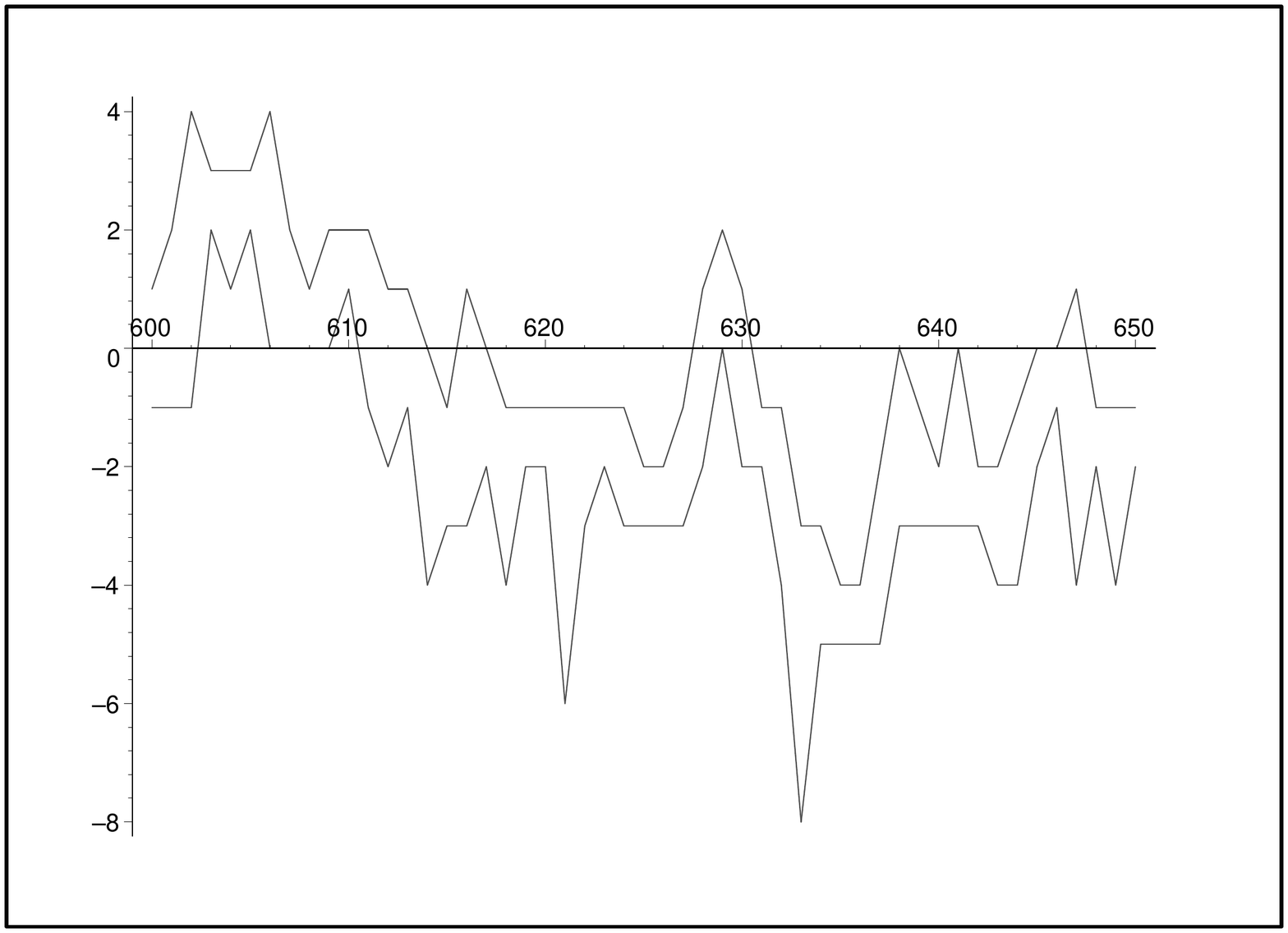}
	\caption{Zoom on our typical polyomino}
	\label{F4}
\end{figure}

Let us illustrate our results by a few figures. We have chosen the cc polyomino as it shows a Markov property (the dcc polyomino is characterized by iid columns).  In Fig.\ref{F1}, we show a simulation of 
$\frac{X(v)-\mu_2 v}{\sqrt{1000}\sig_2},v=1..1000$, In Fig.\ref{F2}, we show a simulation of 
$\frac{Q(v)-\mu_3 v}{\sqrt{1000}\sig_Q},v=1..1000$. The trajectories are strongly oscillating, a classical property of Brownian motions. Fig.\ref{F3} shows a typical polyomino, with its ``filament silhouette".
Fig.\ref{F4} gives a zoom on this polyomino: its width is clearly $\BO(1)$.

\subsection{A comparison with known GF}
Here, we use  Bousquet-M\'elou \cite[(29)]{BM96}\footnote{see Appendix \ref{A3}}. If we set $q=z,x=w$, we obtain
\bals
X&=\frac{wz}{(1-y)(1-yz)}-\frac{w^2z^3(y^2z;z)_2}{(z;z)_1(yz;z)^2_1(yz;z)_2(y^2z;z)_1}+\BO\lp(1-y)^2\rp,\\
W&=-\frac{wz(y^2z;z)_1}{(1-y)(z;z)_1(yz;z)_1}+\BO(1-y).\\
\end{align*} 
The denominator of 
\[\lim_{y=1}\frac{y(1-y)X}{1+W+yX}\]
is given by
\[ z^4(w-1)+z^3(w^2-w+4)-z^2(w+6)+z(w+4)-1\]
which is exactly $h(w,z)$.

 If we set $q=z,y=x=v$,    we obtain
\bals
X&=\frac{vz}{(1-v)(1-vz)}+\sum_{j=2}^\II 
\frac{(-1)^{j+1}v^{j}(1-v)^{2j-4}z^{j(j+1)/2}(v^2z;z)_{2j-2}}{(z;z)_{j-1}(vz;z)_{j-2}(vz;z)^2_{j-1}(vz;z)_j(v^2z;z)_{j-1}},\\
W&=\sum_{j=1}^\II\frac{(-1)^j v^j(1-v)^{2j-3}z^{j(j+1)/2}(v^2;z)_{2j-1}}{(z;z)_j(vz;z)^3_{j-1}(vz;z)_j(v^2;z)_{j-1}}.\\
\end{align*} 
The denominator of 
\[\frac{y(1-y)X}{1+W+yX}\]
gives a function $F(v,z)$ which leads exactly, by Bender's theorems to $\mu_4,\sig_4^2$.
\section{The dc perimeter}
\subsection{The generating functions}
A directed diagonally-convex polyomino (dc) is made of diagonals such that all cells on a diagonal are contiguous and each cell is adjacent to one cell of the previous diagonal. By a rotation of $45^\circ$, this leads to a lattice where the dc is made of contiguous columns such that each cell of each column must be diagonally adjacent to some cell of the previous column.
Note carefully that we don't have here an horizontal perimeter contribution. Note also that this polyomino is different from the  directed and convex polyomino described in Bousquet-M\'elou \cite{BM96}: this last one may have holes in a diagonal. Let us finally remark that our results offer a quite different form from Fereti\'c  and Svrtan  \cite [Thm.5]{FE96} and Fereti\'c  
\cite [Thm.2]{FE02}.

 We extract from \cite{GL297} some relations we had already obtained, starting from $1$ cell: 

\bals
U(k,j)&=(k-j+2)\kl  j\leq k+1\kr,\\
 \Gam(m,\tet,z)&\dq \sum_{j=1} ^\II \tet^j \Fi(m,j,z)= \sum_{j=1} ^\II \tet^j z^j \sum_{k=j-1} ^\II (k-j+2) \Fi(m-1,k,z),\\
 \De(m,\tet,z)&\dq f_1(\tet,z)\De(m-1,1,z)+f_2(\tet,z)\De'(m-1,1,z)+f_3(\tet z)\De(m-1,\tet z,z), m\geq 2,\\
  f_1(\tet,z)& =(-3\tet z+2)/(1-\tet z)^2,f_2(\tet,z)= 1/(1-\tet z), f_3(\tet,z)= (\tet z)^3/(1-\tet z)^2,\De(1,\tet,z)=1,\\ 
	\FI(w,\tet,z)&=\tet [A_1(\xi,\tet,z)+B_{1,1}(\xi,\tet,z) D_1(w,z)+B_{1,2}(\xi,\tet,z)(\tet,z) D_2(w,z)],\mbox{ with }\\
	D_1&=\At_1(\xi,z)+\Bt_{1,1}(\xi,z)D_1+\Bt_{1,2}(\xi,z)D_2,\\
	D_2&=\At_2(\xi,z)+\Bt_{2,1}(\xi,z)D_1+\Bt_{2,2}(\xi,z)D_2,\mbox{ with the following q-analog Bessel functions}\\
	A_1(\xi,\tet,z)&=\xi\sum_{j=0}	^\II \frac{\xi^j\tet^{3j}z^{3j(j+1)/2}}{(\tet z;z)_j^2},\\
	B_{1,1}(\xi,\tet,z)&=\xi\sum_{j=0}	^\II \frac{\xi^j\tet^{3j}z^{3j(j+1)/2}}{(\tet z;z)_j^2}\cdot 
	\frac{-3\tet z^{j+1}+2}{(1-\tet z^{j+1})^2},\\
	B_{1,2}(\xi,\tet,z)&=\xi\sum_{j=0}	^\II \frac{\xi^j\tet^{3j}z^{3j(j+1)/2}}{(\tet z;z)_j^2}\cdot 
	\frac{1}{(1-\tet z^{j+1})},\\
	\At_1(\xi,z)&=\xi\sum_{j=0}	^\II \frac{\xi^j z^{3j(j+1)/2}}{( z;z)_j^2},\\
	\Bt_{1,1}(\xi,z)&=\xi\sum_{j=0}	^\II \frac{\xi^j z^{3j(j+1)/2}}{( z;z)_j^2}\cdot
		\frac{-3 z^{j+1}+2}{(1- z^{j+1})^2},\\ 
		\Bt_{1,2}(\xi,z)&=\xi\sum_{j=0}	^\II \frac{\xi^j z^{3j(j+1)/2}}{( z;z)_j^2}\cdot  \frac{1}{(1- z^{j+1})},\\
		\At_2(\xi,z)&=\xi\sum_{j=0}	^\II \frac{\xi^j z^j z^{3j(j+1)/2} ,}{( z;z)_j^2}[f'_3(z^j,z)A_1(\xi,z^{j+1},z)],\\
		\Bt_{2,1}(\xi,z)&=\xi\sum_{j=0}	^\II \frac{\xi^j  z^j  z^{3j(j+1)/2}}{( z;z)_j^2} 
		[f'_1(z^j,z)+f'_3(z^j,z)B_{1,1}(\xi,z^{j+1},z)] ,\\
		\Bt_{2,2}(\xi,z)&=\xi\sum_{j=0}	^\II \frac{\xi^j  z^j  z^{3j(j+1)/2}}{( z;z)_j^2} 
		[f'_2(z^j,z)+f'_3(z^j,z)B_{1,2}(\xi,z^{j+1},z)] .\\
		&\mbox{This leads to the following expressions}\\
		N_1(w,z)&=-\Bt_{1,2}(w,z)\At_{2}(w,z)-\At_1(w,z)+\Bt_{2,2}(w,z)\At_{1}(w,z),\\
		N_2(w,z)&=\Bt_{1,1}(w,z)\At_{2}(w,z)-\At_2(w,z)-\Bt_{2,1}(w,z)\At_{1}(w,z),\\
		h(w,z)&=\Bt_{1,2}(w,z)\Bt_{2,1}(w,z)-\Bt_{2,2}(w,z)\Bt_{1,1}(w,z)+\Bt_{2,2}(w,z)+\Bt_{1,1}(w,z)-1.\\
		C_2&=C_2(1)=0.3283408377\ldots,\mu_1=0.7660601183\ldots,\mu_2=1.305380578\ldots,\\
		\sig_1^2&=0.1686482431\ldots,\sig_2^2=.3751399028\ldots,\ro=.3756774483\ldots
		\end{align*}
		We can check that $D_1,D_1$ are meromorphic functions for $|w|<1,|z|<1$.  The convergence in the $j$ summations is quite fast. Usually $6$ or $7$ terms are sufficient. To be sure that $\ro$ is the dominant singularity, we can use the principle of the argument of Henrici \cite{HE88}: the number of solutions of an equation $f(z)=0$ that lie inside a simple closed curve $\Gam$, with $f(z)$ analytic inside and on $\Gam$, is equal to the variation of the argument of $f(z)$ along $\Gam$, a quantity also equal to the winding number of the transformed curve $f(\Gam)$ around the origin.  See the application in \cite{GL297}. 
		
		Let us now start  with the case a first column of $i$ cells, (this was not developed in \cite{GL297})
		\bals
		\De(1,\tet,z,i)&=f_0(\tet,z,i)=\tet^{i-1}z^{i-1},\\
		\psi(\xi,\tet,z,i)&=\xi[f_0(\tet,z,i)+f_1(\tet,z) D_1(\xi,z,i)+f_2(\tet,z) D_2(\xi,z),i]+\xi f_3(\tet,z)\psi(\xi,\tet z,z,i),\\ \psi'(\xi,\tet,z,i)&=\xi[f'_0(\tet,z,i)+f'_1(\tet,z) D_1(\xi,z,i)+f'_2(\tet,z) D_2(\xi,z),i]\\
		&+\xi f'_3(\tet,z)\psi(\xi,\tet z,z,i)+ 
		\xi f_3 z(\tet,z)\psi'(\xi,\tet z,z,i),\\
		&\mbox{ we write this as }\\
		\psi(\xi,\tet,z,i)&= \la_1(\tet,z,i)+\mu_{1,1}(\tet,z)\psi(\xi,\tet z,z,i)+\mu_{1,2}(\tet,z)\psi'(\xi,\tet z,z,i),\\
		\psi'(\xi,\tet,z,i)&= \la_2(\tet,z,i)+\mu_{2,1}(\tet,z)\psi(\xi,\tet z,z,i)+\mu_{2,2}(\tet,z)\psi'(\xi,\tet z,z,i),\mbox{ with }\\ 
		\la_1(\tet,z,i)&=\xi[f_0(\tet,z,i)+f_1(\tet,z) D_1(w,z,i)+f_2(\tet,z) D_2(w,z),i],\mu_{1,1}(\tet,z)
		=\xi f_3(\tet,z), \mu_{1,2}(\tet,z)=0,\\ 
		\la_2(\tet,z,i)&=\xi[f'_0(\tet,z,i)+f'_1(\tet,z) D_1(w,z,i)+f'_2(\tet,z) D_2(w,z),i],\\
		\mu_{2,1}(\tet,z)&=
		\xi f'_3(\tet,z), \mu_{2,2}(\tet,z)=\xi f_3(\tet,z)z.\\
		\end{align*}
		Iterating and setting $\tet=1$ leads to
	\bals
				\psi(\xi,\tet,z,i)&=\\
			\la_1(\tet,z,i)&+\mu_{1,1}(\tet,z)\lb \la_1(\sig^{(1)}(\tet),z,i)+\mu_{1,1}(\sig^{(1)}(\tet),z) 
			\lb \la_1(\sig^{(2)}(\tet),z,i) +\mu_{1,1}(\sig^{(2)}(\tet),z)\psi(\xi,\sig^{(3)}(\tet),z,i)\rb\rb,\\
			\sig^{(1)}(\tet)&=\tet z,\sig^{(2)}(\tet)=(\tet z)z=\tet z^2,\sig^{(k)}(\tet)=\tet z^k,\\ 
				D_1(w,z,i)&=\la_1(1,z,i)+\mu_{1,1}(1,z)\la_1(z,z,i)+\mu_{1,1}(1,z)\mu_{1,1}(z,z)\la_1(z^2,z,i)+\ldots\\
				&=\At_1(\xi,z,i)+\Bt_{1,1}(\xi,z)D_1(w,z,i)+\Bt_{1,2}(\xi,z)D_2(w,z,i),\mbox{ with }\\
				A_1(\xi,\tet,z,i)&=\xi\sum_{j=0}	^\II \frac{\xi^j\tet^{3j} z^{3j(j+1)/2}\tet^{i-1}z^{(i-1)(j+1)}}{(\tet z;z)_j^2},\\
				\At_1(\xi,z,i)&=\xi\sum_{j=0}	^\II \frac{\xi^j z^{3j(j+1)/2}z^{(i-1)(j+1)}}{( z;z)_j^2},\\
				&\mbox{ we derive}\\
\psi(\xi,\tet,z,i)&=A_1(\xi,\tet,z,i)+B_{1,1}(\xi,\tet,z)D_1(w,z,i)+B_{1,2}(\xi,\tet,z)D_2(w,z,i),\\
\FI(w,\tet,z,i)&=\tet \psi(wz,\tet,z,i),\mbox{ note that we have no }\xi \mbox{ factor here in front of }\psi\\
 \psi'(\xi,\tet,z,i)&=\\
  \la_2(\tet,z,i)&+\mu_{2,1}(\tet,z)\lb A_1(\xi,\tet z,z,i)+B_{1,1}(\xi,\tet z,z)D_1(w,z,i)+B_{1,2}(\xi,\tet z,z)D_2(w,z,i)\rb \\
			&+\mu_{2,2}(\tet ,z)\psi'(\xi,\tet z,z,i).
	\end{align*}	
	We set $\tet=1$ and rewrite $\psi'(\xi,\tet,z,i)$
\bals	
			D_2(w,z,i)&=H(1,z,i)+H(z,z,i)[\xi zf_3(1,z)]+H(z^2,z,i)[\xi zf_3(1,z)][\xi zf_3(z,z)]+\ldots,\\
\psi'(\xi,\tet,z,i)&=H(\tet,z,i)+\mu_{2,2}(\tet ,z)\psi'(\xi,\tet z,z,i), \mbox{ with }\\
 	H(\tet,z,i)&=\xi\lb f'_0(\tet,z,i)+ f'_1(\tet,z)D_1(w,z,i)+ f'_2(\tet,z)D_2(w,z,i)+f'_3(\tet,z)A_1(\xi,\tet z,z,i)\right.\\
		&\left.+f'_3(\tet,z)B_{1,1}(\xi,\tet z,z)D_1(w,z,i)+f'_3(\tet,z)B_{1,2}(\xi,\tet z,z)D_2(w,z,i)   \rb,\\
		&\mbox{ iterating again }\\
			D_2(w,z,i)&=\At_2(\xi,z,i)+\Bt_{2,1}(\xi,z)D_1(w,z,i)+\Bt_{2,2}(\xi,z)D_2(w,z,i),\mbox{ with }\\
A_2(\xi,\tet,z,i)&=\xi\sum_{j=0}^\II \frac{\xi^j z^j z^{3j(j+1)/2} \tet^{3j}
 [f'_3(\tet z^j,z)A_1(\xi,\tet z^{j+1},z,i)]}{(\tet z;z)_j^2} 
	+\xi\sum_{j=0}	^\II \frac{\xi^j z^j \tet^{3j} z^{3j(j+1)/2}(i-1)\tet^{i-2}z^{(i-2)j}z^{i-1}}{(\tet z;z)_j^2},\\
	\At_2(\xi,z,i)&=\xi\sum_{j=0}	^\II \frac{\xi^j z^j z^{3j(j+1)/2}[f'_3(z^j,z)A_1(\xi, z^{j+1},z,i)]}{( z;z)_j^2}
		+\xi\sum_{j=0}	^\II \frac{\xi^j z^j z^{3j(j+1)/2}(i-1)z^{(i-2)j}z^{i-1}}{( z;z)_j^2},\\
		&\mbox{this gives }\\
 	G(\tet)&=\tet\lb B_{1,1}(\ro,\tet,\ro) N_1(1,\ro)+B_{1,2}(\ro,\tet,\ro) \frac{N_2(1,\ro)}{N_1(1,\ro)} \rb,\mbox{ independent of }i.\\
				\end{align*}	
\subsection{The perimeter conditioned on $n$}

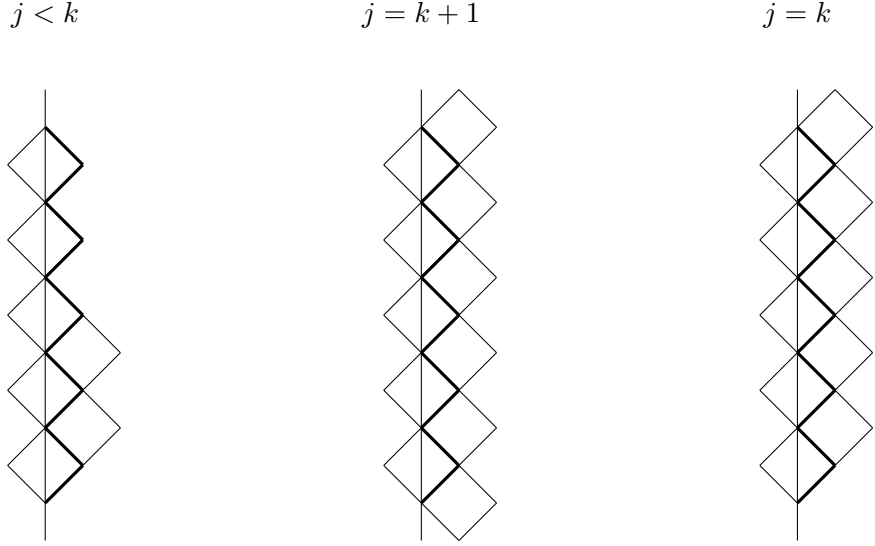
\begin{figure}[htbp]
\begin{center}
\begin{tikzpicture}[set style={{help lines}+=[dashed]}]
  \node at(0,7)[centered]{$j<k$}; \node at(5,7)[centered]{$j=k+1$}; \node at(10,7)[centered]{$j=k$};
  \draw (0,6) --  (0,0); \draw (5,6) --  (5,0); \draw (10,6) --  (10,0);
  \draw (0,0.5) --  (-0.5,1); \draw (-0.5,1) --  (0,1.5);
  \draw (0,1.5) --  (-0.5,2); \draw (-0.5,2) --  (0,2.5);
  \draw (0,2.5) --  (-0.5,3); \draw (-0.5,3) --  (0,3.5);
  \draw (0,3.5) --  (-0.5,4); \draw (-0.5,4) --  (0,4.5);
  \draw (0,4.5) --  (-0.5,5); \draw (-0.5,5) --  (0,5.5);
  \draw[very thick] (0,0.5) --  (0.5,1); \draw[very thick] (0.5,1) --  (0,1.5);
  \draw[very thick] (0,1.5) --  (0.5,2); \draw[very thick] (0.5,2) --  (0,2.5);
  \draw[very thick] (0,2.5) --  (0.5,3); \draw[very thick] (0.5,3) --  (0,3.5);
  \draw[very thick] (0,3.5) --  (0.5,4); \draw[very thick] (0.5,4) --  (0,4.5);
  \draw[very thick] (0,4.5) --  (0.5,5); \draw[very thick] (0.5,5) --  (0,5.5);
  \draw (0.5,1) --  (1,1.5); \draw (1,1.5) --  (0.5,2);
  \draw (0.5,2) --  (1,2.5); \draw (1,2.5) --  (0.5,3);
  \draw (5,0.5) --  (4.5,1); \draw (4.5,1) --  (5,1.5);
  \draw (5,1.5) --  (4.5,2); \draw (4.5,2) --  (5,2.5);
  \draw (5,2.5) --  (4.5,3); \draw (4.5,3) --  (5,3.5);
  \draw (5,3.5) --  (4.5,4); \draw (4.5,4) --  (5,4.5);
  \draw (5,4.5) --  (4.5,5); \draw (4.5,5) --  (5,5.5);
  \draw[very thick] (5,0.5) --  (5.5,1); \draw[very thick] (5.5,1) --  (5,1.5);
  \draw[very thick] (5,1.5) --  (5.5,2); \draw[very thick] (5.5,2) --  (5,2.5);
  \draw[very thick] (5,2.5) --  (5.5,3); \draw[very thick] (5.5,3) --  (5,3.5);
  \draw[very thick] (5,3.5) --  (5.5,4); \draw[very thick] (5.5,4) --  (5,4.5);
  \draw[very thick] (5,4.5) --  (5.5,5); \draw[very thick] (5.5,5) --  (5,5.5);
  \draw (5.5,0) --  (6,0.5); \draw (6,0.5) --  (5.5,1);
  \draw (5.5,1) --  (6,1.5); \draw (6,1.5) --  (5.5,2);
  \draw (5.5,2) --  (6,2.5); \draw (6,2.5) --  (5.5,3);
  \draw (5.5,3) --  (6,3.5); \draw (6,3.5) --  (5.5,4);
  \draw (5.5,4) --  (6,4.5); \draw (6,4.5) --  (5.5,5);
  \draw (5.5,5) --  (6,5.5); \draw (6,5.5) --  (5.5,6);
  \draw (5,0.5) --  (5.5,0); \draw (5,5.5) --  (5.5,6);
  \draw (10,0.5) --  (9.5,1); \draw (9.5,1) --  (10,1.5);
  \draw (10,1.5) --  (9.5,2); \draw (9.5,2) --  (10,2.5);
  \draw (10,2.5) --  (9.5,3); \draw (9.5,3) --  (10,3.5);
  \draw (10,3.5) --  (9.5,4); \draw (9.5,4) --  (10,4.5);
  \draw (10,4.5) --  (9.5,5); \draw (9.5,5) --  (10,5.5);
  \draw[very thick] (10,0.5) --  (10.5,1); \draw[very thick] (10.5,1) --  (10,1.5);
  \draw[very thick] (10,1.5) --  (10.5,2); \draw[very thick] (10.5,2) --  (10,2.5);
  \draw[very thick] (10,2.5) --  (10.5,3); \draw[very thick] (10.5,3) --  (10,3.5);
  \draw[very thick] (10,3.5) --  (10.5,4); \draw[very thick] (10.5,4) --  (10,4.5);
  \draw[very thick] (10,4.5) --  (10.5,5); \draw[very thick] (10.5,5) --  (10,5.5);
  \draw (10.5,1) --  (11,1.5); \draw (11,1.5) --  (10.5,2);
  \draw (10.5,2) --  (11,2.5); \draw (11,2.5) --  (10.5,3);
  \draw (10.5,3) --  (11,3.5); \draw (11,3.5) --  (10.5,4);
  \draw (10.5,4) --  (11,4.5); \draw (11,4.5) --  (10.5,5);
  \draw (10.5,5) --  (11,5.5); \draw (11,5.5) --  (10.5,6);
  \draw (10,5.5) --  (10.5,6);
 \end{tikzpicture}
\end{center}
\caption{Three possibilities  of a dc polyomino .}
\label{F7}
\end{figure}

 For dc, we have the following 
possibilities: see Fig.\ref{F7},  

\bals
&\mbox{ The possibility function is given by }U(k,j)=(k-j+2)\kl  j\leq k+1\kr,\\
&\mbox{ we note that, if we add a cell at the upper or lower  previous cell position of column }d-1,\\
& \mbox{ we increase the perimeter by }2,
\mbox{ otherwize, the perimeter doesn't change.  }\\
 &\mbox{Hence the first moments are computed as}\\
 E_w(k,j)&=\frac{1}{k-j+2}(1\cdot 2+1\cdot 2)=\frac{4}{k-j+2},\\
\mbox{ if }k=j+1,&\B E_{w^2}(k,j)=16,\\
\mbox{ if }k\leq j,&\B E_{w^2}(k,j)=\frac{1}{k-j+2}(1\cdot 4+1\cdot 4)=\frac{8}{k-j+2}.\\
&\mbox{Finally, we apply the results from previous sections},\\
&\mbox{we make the following substitutions (note that we divide again by }w)\\
		\mbox{in }	\Xi_3(F_1,F_2)&:M_1(i,\ell):=\mbox{ subs }[\xi=zw,w=1,z=\ro],\mbox{ in }[\tet^\ell]\tet A_1(\xi,\tet,z,i)/w,\\
	\mbox{in }		\Xi_4(F_1,F_2)&:\Fi(w,\tet,\ro,i):=
\mbox{ subs }[\xi=\ro w ]\mbox{ in }\tet [B_{1,1}(\xi,\tet,\ro) N_1(w,\ro,i) 	+B_{1,2}(\xi,\tet,\ro) N_2(w,\ro,i)]/w,\\
&\mbox{and we obtain }\\
\mu_3&=2.2705856475\ldots,\ro_{X,Q}=0.5713021769\ldots,\mu_4=1.7394051099\ldots,\sig_3^2=.9808725500\ldots\\
\sig_4^2&=0.38150889574\ldots,\sig_Q^2=0.362055589\ldots,\mbox{ again } \sig_X^2 \equiv \sig_2^2.\\
\end{align*}
Again, we have made extensive simulations to check our results. The fit is quite good.
\section{The staircase perimeter}
A staircase (or parallelogram) polyomino (st), is made of contiguous columns such that the base cell of each column  must be adjacent to some cell of the previous column and the top cell of each column must be adjacent to some cell of the next column. 
\subsection{The generating functions}

We proceed as in the previous sections. Starting with a first column of size $i$, we have
\bals
U(k,j)&= j \mbox{ if }j\leq k,U(k,j)=k \mbox{ if }j> k.\\
&\mbox{The preliminary relations are}\\
\Gam(m,\tet,z)&=\sum_{j=1} ^\II \tet^j \Fi(m,j,z)\dq \sum_{j=1} ^\II \tet^j z^j\lb \sum_{k=1} ^{j-1} k \Fi(m-1,k,z)
+ \sum_{k=j} ^\II  j \Fi(m-1,k,z)\rb,\\
\De(m,\tet,z)&=f_1(\tet,z)\De(m-1,1,z)+f_2(\tet,z)\De(m-1,\tet z,z),\\
f_1(\tet,z)&=\frac{1}{(1-\tet z)^2},f_2(\tet,z)=-\frac{\tet z}{(1-\tet z)^2},\De(1,\tet,z)=f_0(\tet,z)=\tet^{i-1}z^{i-1},\\
\FI(w,\tet,z,i)&=\tet \psi(wz,\tet,z,i)
=\tet[A_1(w,\tet,z,i)+B_{1}(w,\tet,z)D_1(w,z,i)],\mbox{ with }\\
A_1(\xi,\tet,z,i)&=\xi\sum_{j=0}	^\II \frac{\xi^j\tet^j z^{3j(j+1)/2}(-1)^j}{(\tet z;z)_j^2}\tet^{i-1}z^{(i-1)(j+1)},\\
\At_1(\xi,z,i)&=\xi\sum_{j=0}	^\II \frac{\xi^j z^{3j(j+1)/2}(-1)^j}{( z;z)_j^2}z^{(i-1)(j+1)},\\
B_1(\xi,\tet,z)&=\xi\sum_{j=0}	^\II \frac{\xi^j\tet^j z^{3j(j+1)/2}(-1)^j}{(\tet z;z)_j^2}\frac{1}{(1-\tet z^{j+1})^2},\\
\Bt_1(\xi,z)&=\xi\sum_{j=0}	^\II \frac{\xi^j z^{3j(j+1)/2}(-1)^j}{( z;z)_j^2}\frac{1}{(1- z^{j+1})^2},\\
D_1(\xi,z,i)&=\frac{\At_1(\xi,z,i)}{1-\Bt_1(\xi,z)}=\frac{N_1(w,z,i)}{h(w,z)},C_2(j)=-\frac{N_1(1,\ro,j)}{h_w(1,\ro)}.\\
&\mbox{We obtain }\\
G(\tet)&=\frac{\tet B_{1}(1,\tet,\ro)}{\Bt_{1}(1,\ro)},\mbox{ but }\Bt_{1}(1,\ro)=1 \mbox{ by }h(1,\ro)=0,\\
\mu_1&=0.4208810078\ldots,\sig_1^2=0.2080626954\ldots,C_2=0.3060622477\ldots=C_2(1),\ro = .4330619231\ldots,\\
\mu_2&=2.3759684098\ldots,\sig_2^2=2.7907198037\ldots
\end{align*}
\subsection{The Markov chain}
\bals
\mbox{if }& j\leq k, \mbox{ we denote   any function } F(k,j) \mbox{ by } F^m(k,j)\mbox{ and by }  F^p(k,j)\mbox{ otherwise, this gives },\\
\Pi^m(k,j)&\dq \pi(k) j C_2(j)/\pi_2(k),j\leq k,\\
\Pi^p(k,j)&\dq \pi(k) k C_2(j)/\pi_2(k),j>k,\\
\pi_2(k)&=\pi(k)\lb \sum_{j=1}^k j C_2(j)+\sum_{j=k+1}^\II k C_2(j)\rb,\\
\pi_2(j)&=C_2(j)\lb \sum_{k=j}^\II \pi(k) j +\sum_{k=1}^{j-1} \pi(k) k \rb.\\
 \end{align*}
\subsection{The perimeter conditioned on $n$}

\begin{figure}[htbp]
\begin{center}
\begin{tikzpicture}[set style={{help lines}+=[dashed]}]
 \draw (2,0) rectangle (1,5);
 \draw (1,4) --  (2,4);
 \node[circle,fill=black!100,inner sep=0.05cm](s1)at(1.5,4.5){};
 \draw (1,2) --  (2,2);
 \node[circle,fill=black!100,inner sep=0.05cm](s1)at(1.5,1.5){}; 
 \draw (1,1) --  (2,1);
 \node[circle,fill=black!100,inner sep=0.05cm](s1)at(1.5,0.5){};
 \draw (-2,-4) rectangle (-1,2);
 \node[circle,fill=black!100,inner sep=0.05cm](s1)at(-1.5,1.5){};
 \draw (-2,1) --  (-1,1);
 \draw (-2,0) --  (-1,0);
 \node[circle,fill=black!100,inner sep=0.05cm](s1)at(-1.5,0.5){};
  \draw[dashed] (-2,5) --  (1,5);
  \draw[dashed] (-1,2) --  (2,2);
  \draw[dashed] (-2,0) --  (2,0);
  \draw[dashed] (-2,-4) --  (2,-4);
\draw [decorate,decoration={brace,amplitude=10pt},xshift=-4pt,yshift=0pt]
(-2,2) -- (-2,5) node [black,midway,xshift=-0.6cm] 
{\footnotesize $w$};
\draw [decorate,decoration={brace,amplitude=10pt,mirror,raise=4pt},yshift=0pt]
(2,-4) -- (2,0) node [black,midway,xshift=0.8cm] {\footnotesize
$z$};
\end{tikzpicture}
\end{center}
\caption{Two columns of a st polyomino and their related parameters.}
\label{F8}
\end{figure}
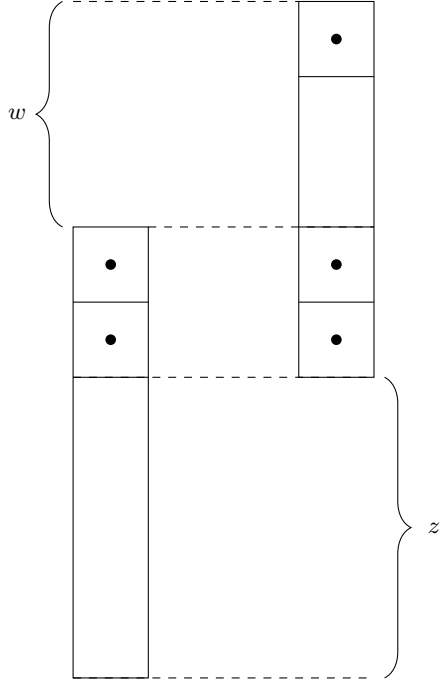

 For st,   we have the following 
notations and probabilistic relations: see Fig.\ref{F8},  

\bals
  w_d&:=\mbox{ upper  cell position of column }d-\mbox{ upper  cell position of column }d-1,\\
		z_d&:=\mbox{ lower  cell position of column }d-\mbox{ lower  cell position of column }d-1,\\
		\mbox{if }& j\leq k, w_d \geq 0\mbox { takes uniformly values on }(0,j-1),z_d \geq 0\mbox { takes uniformly values on }(k-j,k-1),\\
		\mbox{if }& j> k, w_d \geq 0\mbox { takes uniformly values on }(j-k,j-1),z_d \geq 0\mbox { takes uniformly values on }(0,k-1)\\
		T_d&:=w_d+z_d,w=j-k+z.\\
		&\mbox{The necessary moments are computed as follows }\\
		E^m_z(k,j)&=1/j \sum_{u=1}^{j}(k-u)=k-1/2 j-1/2,E^p_z(k,j)=1/k \sum_{u=0}^{k-1}u=1/2 k-1/2,\\
		E^m_w(k,j)&=1/j \sum_{u=0}^{j-1}u=1/2 j-1/2, E^p_w(k,j)=1/k \sum_{u=1}^{k}(j-u)=j-1/2 k-1/2,\\
		E^m_{z^2}(k,j)&=1/j \sum_{u=1}^{j}(k-u)^2=k^2-k j-k+1/3 j^2+1/2 j+1/6,\\
		E^m_{w^2}(k,j)&=1/j  \sum_{u=0}^{j-1}u^2=1/3 j^2-1/2 j+1/6,\\
		E^p_{z^2}(k,j)&=1/k \sum_{u=0}^{k-1}u^2=1/3 k^2-1/2 k+1/6,\\
		E^p_{w^2}(k,j)&=1/k \sum_{u=1}^{k}(j-u)^2=j^2-k j-j+1/3 k^2+1/2 k+1/6,\\
			E^m_{zw}(k,j)&=1/j \sum_{u=1}^{j} (k-u) (j-u) =1/2 k j-1/2 k-1/6 j^2+1/6,\\
			E^p_{zw}(k,j)&=1/k \sum_{u=0}^{k-1} u(j-k+u)=1/2 k j-1/2 j-1/6 k^2+1/6,\\
			E^m_{zaw}(k,j)&=E^m_z(k,j)+E^m_w(k,j)=k-1,E^p_{zaw}(k,j)=E^p_z(k,j)+E^p_w(k,j)=j-1,\\
			E^m_{z^2aw^2}(k,j)&=E^m_{z^2}(k,j)+E^m_{w^2}(k,j)=k^2-k j-k+2/3 j^2+1/3,\\
			E^p_{z^2aw^2}(k,j)&=E^p_{z^2}(k,j)+E^p_{w^2}(k,j)=2/3 k^2+1/3+j^2-k j-j,\\
			&\mbox{again  we have some symmetry here: }E^m_{.}(k,j)\equiv E^p_{.}(j,k),\\
			&\mbox{ for } y^m \mbox{ a random variable depending on } k,j, \mbox{ with mean }E^m_y(k,j),\\
			&\mbox{ we set }\Eb^m(y)\dq \sum_{u=1}^\II \sum_{j=1}^u \pi(u) j C_2(j) E^m_y(k,j),\\
			&\mbox{ for } y^p \mbox{ a random variable depending on } k,j, \mbox{ with mean }E^p_y(k,j),\\
			&\mbox{ we set }\Eb^p(y)\dq \sum_{u=1}^\II \sum_{j=u+1}^u \pi(u) u C_2(j)E^m_y(k,j),\\
			\Eb(y)&=\Eb^m(y)+\Eb^p(y),\mbox{ and similarly for functions }F^m(k,j),F^p(k,j),\\
			\sig_3^2&=\Eb\lp z^2+w^2+2 zw\rp-\mu_3^2=3.3102701914\ldots,\mu_3=\Eb(z+w)=2.
			\end{align*}
			We use the previous relations as follows: we make the following substitutions 
			\bals
		\mbox{in }	\Xi_3(F_1,F_2)&:M_1(i,\ell):=\mbox{ subs }[\xi=zw,w=1,z=\ro],\mbox{ in }[\tet^\ell]\tet A_1(\xi,\tet,z,i)/w,\\
			\mbox{in }	\Xi_4(F_1,F_2)&:\Fi(w,\tet,\ro,i):=\mbox{ subs }[\xi=\ro w ]\mbox{ in }\tet [B_{1}(\xi,\tet,\ro) N_1(w,\ro,i)]/w,\\
				\Xi_5(F_1,F_2)&=	\left\{ \sum_{u=1}^\II \pi(u)\sum_{j=1}^u  j F_1^m(u,j)+\sum_{u=1}^\II \pi(u)\sum_{j=u+1}^\II  u F_1^p(u,j)\right\} \cdot \\
				&\cdot \left\{ \sum_{\ell=1}^\II [M_1(j,\ell)+M_2(j,\ell)]\sum_{k=1}^\ell k C_2(k)  F_2^m(\ell,k) +\sum_{\ell=1}^\II[M_1(j,\ell)+M_2(j,\ell)] \sum_{k=\ell+1}^\II \ell C_2(k)  F_2^p(\ell,k)\right\},\\
				\sig_Q^2&=\sig_3^2+2 \Xi_5(E_{zaw},E_{zaw})=6.199368675211\ldots,\\
				\sig_X^2&=\sig_x^2+2\Xi_5(j,k)\equiv \sig_2^2,\\
				C_{X,Q}&=\Eb(j\cdot E_{zaw})-\mu_2\mu_3+\lb \Xi_5(j,E_{zaw})+\Xi_5(E_{zaw},k) \rb,\mu_4=1.683524031\ldots,\\
				\sig_4^2&=0.7198047885\ldots,\ro_{X,Q}=0.8853121502\ldots,\sig_3^2=3.3102701914\ldots
\end{align*} 
As previously, we have made extensive simulations to check our results. The fit is quite good.
\subsection{A comparison with known GF}
Now we use  Bousquet-M\'elou \cite[Thm. 3.2]{BM96}. If we set $q=z,y=1,x=w$ 
in the denominator $J_0(1)$, we obtain a function
\[\sum_{n=0}^\II \frac{(-1)^n w^n z^{n(n+1)/2}}{(z;z)^2_n}\]
which is another form of
\[h(w,z)=1-\Bt_1(\xi,z).\]
 If we set $q=z,y=x=v$,    we obtain
\[F(v,z)=\sum_{n=0}^\II \frac{(-1)^n v^n z^{n(n+1)/2}}{(z;z)_n (vz;z)_n}\]
which leads exactly, by Bender's theorems to $\mu_4,\sig_4^2$. This last analysis is also given in 
Flajolet and Sedgewick \cite [Prop.IX.11]{FLSe09}
\section{The escalier perimeter}
The escalier polyomino (es) is made of contiguous columns with all base cells at the same level, such that, if the size of a column is $k$ and the size of the next column is $j$, we must have $j\geq k-1$. Let us  remark that our results offer a quite different form from   Fereti\'c \cite [Prop.2]{FE02}.
\subsection{The generating functions and Markov chain}

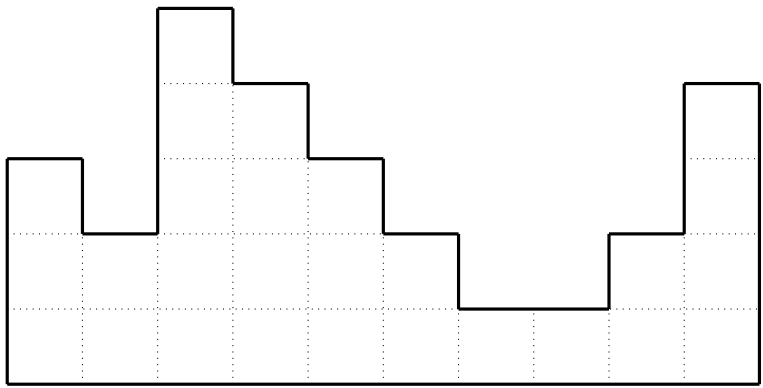
\begin{figure}[htbp]
\begin{center}
\begin{tikzpicture}[set style={{help lines}+=[dashed]}]
 \draw[very thick] (0,0) --  (10,0);
 \draw[very thick]  (0,0) --  (0,3);
 \draw[very thick]  (10,0) --  (10,4);
 \draw[very thick] (0,3) --  (1,3); \draw[very thick] (1,3) --  (1,2);
 \draw[very thick] (1,2) --  (2,2); \draw[very thick] (2,2) --  (2,5);
 \draw[very thick] (2,5) --  (3,5); \draw[very thick] (3,5) --  (3,4);
 \draw[very thick] (3,4) --  (4,4);
 \draw[very thick] (4,4) --  (4,3); \draw[very thick] (4,3) --  (5,3);
 \draw[very thick] (5,3) --  (5,2); \draw[very thick] (5,2) --  (6,2);
 \draw[very thick] (6,2) --  (6,1); \draw[very thick] (6,1) --  (8,1);
 \draw[very thick] (8,1) --  (8,2); \draw[very thick] (8,2) --  (9,2);
 \draw[very thick] (9,2) --  (9,4); \draw[very thick] (9,4) --  (10,4);
  \draw[dotted] (0,1) --  (10,1);
  \draw[dotted] (0,2) --  (1,2);  \draw[dotted] (2,2) --  (5,2);  \draw[dotted] (9,2) --  (10,2);
  \draw[dotted] (2,3) --  (4,3);  \draw[dotted] (9,3) --  (10,3); \draw[dotted] (2,4) --  (3,4);
  \draw[dotted] (1,0) --  (1,2);  \draw[dotted] (2,0) --  (2,2); \draw[dotted] (3,0) --  (3,4);
  \draw[dotted] (4,0) --  (4,3);  \draw[dotted] (5,0) --  (5,2); \draw[dotted] (6,0) --  (6,1);
  \draw[dotted] (7,0) --  (7,1); \draw[dotted] (8,0) --  (8,1); \draw[dotted] (9,0) --  (9,2);
 \end{tikzpicture}
\end{center}
\caption{A typical es polyomino.}
\label{F9}
\end{figure}

 For es, we have the following 
typical example: see Fig.\ref{F9},  

We have here 
\bals
U(k,j)&=\kl j\geq k-1 \kr.\\
&\mbox{ Starting with a first column of size }i,\\
\Fi(m,j,z)&=z^j\sum_{k=1}^{j+1}\Fi(m-1,k,z),\Fi(1,i,z)=z^i,\\ 
\Gam(m,\tet,z)&=f_1(\tet,z)\Gam(m-1,\tet z,z)-\Fi(m-1,1,z),\mbox{ with }f_1(\tet,z)=\frac{1}{\tet z(1-\tet z)},
\Gam(1,\tet,z)=\tet^i z^i.\\ 
\end{align*}
This polyomino is quite different from other polyominoes and our usual techniques do not work anymore. The presence  of $\tet$
 in the denominator of $f_1(\tet,z)$ excludes a direct iteration procedure. We must turn to another approach. We are indebted to H.Prodinger and S.Wagner for providing a new analysis: \cite{PW18}. First of all, let us change  the notations. 
Starting with a first column of size $1$, we have
\bals
f(n, j,z) &=\Fi (n + 1, j + 1,z),\\
f(0, 0,z)& = z, f(0, k,z) = 0 \mbox{ for } k \geq 1,\\
\Fi(m, j,z)& = z^j\sum_{1\leq k\leq j+1}\Fi(m - 1, k,z),\\
\Fi(m + 1, j,z) &= z^j \sum_{0\leq k\leq j} \Fi(m, k + 1,z),\\
\Fi(m + 1, j + 1,z)& = z^{j+1} \sum_{0\leq k\leq j+1}\Fi(m, k + 1,z),\\
f(n, j,z)& = z^{j+1}\sum_{0\leq k\leq j+1}f(n-1  , k,z).\\
&\mbox{Now set}\\
F(x, y,z)& =\sum_{n\geq 0} \sum_{j\geq 0}x^n y^j f(n, j,z),\\
\FI(w,\tet,z)&=w\tet F(w,\tet,z),\\
F(x, y,z)& = z +\sum_{n\geq 1} \sum_{j\geq 0}x^n y^j z^{j+1} \sum_{0\leq k\leq j+1}f(n - 1, k,z),\\
&= z + x \sum_{n\geq 0}x^n \sum_{j\geq 0}y^j z^{j+1} \sum_{0\leq k\leq j+1}f(n, k,z)\\
&= z + x \sum_{n\geq 0}x^n  \sum_{j\geq 0}y^j z^{j+1}f(n, 0,z)\\
&+ x\sum_{n\geq 0} x^n \sum_{j\geq 0}y^j z^{j+1}\sum_{1\leq k\leq j+1}f(n, k,z)\\
&= z +\frac{zx}{1 - zy}\sum_{n\geq 0}x^nf(n, 0,z)\\
&+ x\sum_{n\geq 0}x^n\sum_{k\geq 1}f(n, k,z)\sum_{j\geq k-1}y^j z^{j+1}\\
&= z +\frac{zx}{1 - zy}F(x, 0,z)+
\frac{x}{y(1 - zy)}\sum_{n\geq 0}x^n\sum_{k\geq 1}f(n, k,z)(zy)^k\\
&= z +\frac{zx}{1 - zy}F(x, 0,z)+\frac{x}{y(1 - zy)}\sum_{n\geq 0}x^n\sum_{k\geq 0}f(n, k,z)(zy)^k
- \frac{x}{y(1 - zy)}\sum_{n\geq 0}x^n f(n, 0,z)\\
&= z +\frac{zx}{1 - zy}F(x, 0,z)+
\frac{x}{y(1 - zy)}F(x, zy,z) -\frac{x}{y(1 - zy)}F(x, 0,z)\\
&= z -\frac{x}{y}F(x, 0,z) +\frac{x}{y(1 - zy)}F(x, zy,z)\\
F(x, y,z)& =z -\frac{x}{y}F(x, 0,z) +\frac{x}{y(1 - zy)}F(x, zy,z).\\
\end{align*}

Let $F(x,y,z)$ be the unique function that is analytic in $x$ and $y$ around $(0,0)$ and satisfies
$$F(x,y,z) = z - \frac{x}{y} F(x,0,z) + \frac{x}{y(1-zy)} F(x,zy,z).$$
It turns out that we have
$$F(x,y,z) = z \Big (1 + \sum_{n \geq 0} z^{-n(n+1)/2} x^{-n} y^n (Q_n K - P_n) \Big),$$
where $K$ is the continued fraction
$$K = \cfrac{1}{1-\cfrac{zx}{1-\cfrac{z^2x}{1-\cfrac{z^3x}{\ddots}}}} = \frac{\sum_{j \geq 0} \frac{(-1)^j z^{j^2+j} x^j}{(z;z)_j}}{\sum_{j\geq 0} \frac{(-1)^j z^{j^2} x^j}{(z;z)_j}},$$
and $P_n$ and $Q_n$ are numerators and denominators of the convergents of $K$:
$$\cfrac{1}{1-\cfrac{zx}{1-\cfrac{z^2x}{1-\cfrac{z^3x}{\cfrac{\ddots}{1-\cfrac{z^{n-1}x}{1-z^n x}}}}}} = \frac{P_n}{Q_n}.$$
The recursions
$$P_n = P_{n-1} - z^nx P_{n-2}$$
and
$$Q_n = Q_{n-1} - z^nx Q_{n-2}$$
hold with initial values $P_{-1} = 0$ and $P_0 = Q_{-1} = Q_0 = 1$, and we have the explicit formulas
$$P_n = \sum_{j \geq 0} \QB{n-j}{j} (-1)^j z^{j^2+j} x^j$$
as well as
$$Q_n = \sum_{j \geq 0} \QB{n+1-j}{j} (-1)^j z^{j^2} x^j.$$
We also have the recursion
$$Q_n K - P_n = (Q_{n-1}K-P_{n-1}) \cdot \cfrac{z^{n+1} x}{1-\cfrac{z^{n+2}x}{1-\cfrac{z^{n+3}x}{\ddots}}}\,,$$
which gives us
$$Q_n K - P_n = K \cdot \prod_{k=1}^{n+1} \cfrac{z^k x}{1-\cfrac{z^{k+1}x}{1-\cfrac{z^{k+2}x}{\ddots}}}\,.$$
It follows that
$$z^{-n(n+1)/2} x^{-n} (Q_n K - P_n) = K x z^{n+1} \cdot \prod_{k=1}^{n+1} \cfrac{1}{1-\cfrac{z^{k+1} x}{1-\cfrac{z^{k+2}x}{\ddots}}}\,.$$
With
$$P = \sum_{j \geq 0} \frac{(-1)^j z^{j^2+j} x^j}{(z;z)_j}$$
and
$$Q(x,z) = \sum_{j \geq 0} \frac{(-1)^j z^{j^2} x^j}{(z;z)_j}$$
(so that $K = P/Q$), we also have
$$Q_n P - P_n Q =\sum_{k \geq 0}\frac{(-1)^k}{(z;z)_k}z^{(k+1)(k+n+1)+n(n+1)/2}x^{k+n+1},$$
which ultimately yields
$$F(x,y,z) = z + \frac{z}{Q(x,z)} \sum_{n \geq 0} y^n \sum_{k \geq 0}\frac{(-1)^k}{(z;z)_k}z^{(k+1)(k+n+1)}x^{k+1}.$$
and
\[h(w,z)=Q(w,z).\]
In particular,
$$[y^n] F(x,y,z) = \frac{z \sum_{k \geq 0}\frac{(-1)^k}{(z;z)_k}z^{(k+1)(k+n+1)}x^{k+1}}{\sum_{k \geq 0} \frac{(-1)^k z^{k^2} x^k}{(z;z)_k}}$$
for $n \geq 1$.

It is worthwhile to consider the limit as $z \to 1$. In this case, the original functional equation becomes
$$F(x,y,1) = 1 - \frac{x}{y} F(x,0,1) + \frac{x}{y(1-y)} F(x,y,1),$$
whose solution is
$$F(x,y,1) = \frac{1- \frac{xF(x,0,1)}{y}}{1-\frac{x}{y(1-y)}}.$$
Taking $y \to 0$ only yields the trivial identity $F(x,0,1) = F(x,0,1)$ here, but we note that
$$F(x,0,1) = \lim_{z \to 1} zK = \lim_{z \to 1} \cfrac{z}{1-\cfrac{zx}{1-\cfrac{z^2x}{1-\cfrac{z^3x}{\ddots}}}} = \cfrac{1}{1-\cfrac{x}{1-\cfrac{x}{1-\cfrac{x}{\ddots}}}} = \frac{1-\sqrt{1-4x}}{2x}$$
in this case, which yields
$$F(x,y,1) = \frac{(1-y)(1-2y-\sqrt{1-4x})}{2(x-y+y^2)}.$$
The coefficients are given by the generalised Catalan numbers
$$[x^n y^m] F(x,y,1) = \frac{(m+2)(2n+m-1)!}{(n-1)!(n+m+1)!}.$$

With the more general initial condition (originally stated as $\varphi(1,i,z) = z^i$), we obtain the analogous functional equation
\bals
F(x,y,z,i)& = y^{r}z^{i} - \frac{x}{y} F(x,0,z,i) + \frac{x}{y(1-zy)} F(x,zy,z,i),\\
\FI(w,\tet,z,i)&=w\tet F(w,\tet,z,i),\\
\end{align*}
Set $r = i-1$. The solution to this equation is now
$$F(x,y,z,i) = z^{r+1} \Big( y^r + \sum_{n=0}^{r-1} z^{r(r+1)/2-n(n+1)/2} x^{r-n} y^n \frac{Q_n(x,z)Q(z^{r+1}x,z)}{Q(x,z)} + \sum_{n=r}^{\infty} z^{n+1} x y^n \frac{Q_{r-1}(x,z)Q(z^{n+2}x,z)}{Q(x,z)} \Big).$$
So we find that
$$[y^n] F(x,y,z,i) = \begin{cases} z^{(r+1)(r+2)/2-n(n+1)/2}x^{r-n} \frac{Q_n(x,z)Q(z^{r+1}x,z)}{Q(x,z)} & n < r, \\ z^{n+r+2} x \frac{Q_{r-1}(x,z)Q(z^{n+2}x,z)}{Q(x,z)} + \kl n=r \kr z^{r+1} & n \geq r.\end{cases}$$
\bals
&\mbox{for the case }i=1, \mbox{ we have},\\
\FI(w,\tet,z)&=w\tet z+\frac{1}{Q(w,z)}w\tet z \sum_{n \geq 0} \tet^n H_n(w,z),h(w,z)=Q(w,z)\\
H_n(w,z)&:=\sum_{k \geq 0}\frac{(-1)^k}{(z;z)_k}z^{(k+1)(k+n+1)}w^{k+1},\\
S(w,z)&:=wz\sum_{n \geq 0}  H_n(w,z)=wz\sum_{k \geq 0}\frac{(-1)^k}{(z;z)_k}\frac{z^{(k+1)^2}}{1-z^{k+1}}w^{k+1},
C_2=-\frac{S(1,\ro)}{h_w(1,\ro)}\\
G(\tet)&=\frac{\tet \ro \sum_{n \geq 0} \tet^n H_n(1,\ro)}{\ro \sum_{n \geq 0}  H_n(1,\ro)},
\sum_{n \geq 0}  H_n(1,\ro)=1\mbox{ as } h(1,\ro)=0.\\
&\mbox{for the general case }i\geq 1, \mbox{ we have},\\
\FI(w,\tet,z,i)&=w\tet^i z^i+\frac{1}{Q(w,z)}w\tet z^i\lb 
\sum_{n=0}^{i-2} z^{i(i-1)/2-n(n+1)/2} w^{i-1-n} \tet^n Q_n(w,z)Q(z^{i}w,z)\right.\\
&\left.  + \sum_{n=i-1}^{\infty} z^{n+1} w \tet^n Q_{i-2}(w,z)Q(z^{n+2}w,z) \rb,\\
S(w,z,i)&:=w z^i\lb 
\sum_{n=0}^{i-2} z^{i(i-1)/2-n(n+1)/2} w^{i-1-n}  Q_n(w,z)Q(z^{i}w,z)\right.\\
&\left.  + \sum_{n=i-1}^{\infty} z^{n+1} w  Q_{i-2}(w,z)Q(z^{n+2}w,z) \rb,C_2(j)=-\frac{S(1,\ro,j)}{h_w(1,\ro)},\\
P(k,j)&=\pi(k)U(k,j)C_2(j),U(k,j)= \kl k\leq j+1 \kr,\\
\pi_2(k)&=\pi(k)\sum_{j=k-1}^\II C_2(j),k>1,\pi_2(1)=\pi(1)\sum_{j=1}^\II C_2(j),\\
\Pi(1,j)&=\pi(1)C_2(j)/\pi_2(1),j\geq 1,\Pi(k,j)=\pi(k)C_2(j)/\pi_2(k),j\geq 1,k\geq 2,j\geq k-1,\\
\sig_1^2&=0.2290348188\ldots,\mu_1=0.6149126319\ldots,\ro=0.5761487691\ldots,\mu_2=1.626247287\ldots\\
C_2&=C_2(1)=0.8600102250\ldots,\sig_2^2=.9850567845\ldots,\sig_3^2=.3631554767\ldots
\end{align*}
\subsection{The perimeter conditioned on $n$}
We have here 
\bals
U(k,j)&=\kl j\geq k-1 \kr,w=j-k\mbox{ (we have no } z \mbox{ here}),T_d=|w_d|.\\
&\mbox{ We use the notations:}\\
j=k-1&: E^m_w(k,j)=1,E^m_{w^2}(k,j)=1,j\geq k:E^p_w(k,j)=j-k,E^p_{w^2}(k,j)=(j-k)^2,\\
\Eb^m(y)&=\sum_{k=2}^\II \pi(k)C_2(k-1)E_y^m(k,k-1),\\
\Eb^p(y)&=\sum_{k=1}^\II \pi(k)\sum_{j=k}^\II C_2(j)E_y^p(k,j),\\
\Eb(y)&=\Eb^m(y)+\Eb^p(y),\\
\mu_3&=\Eb(w),\sig_3^2=\Eb(w^2)-\mu_3^2.\\
&\mbox{We use previous relations as follows }\\
\Xi_3(F_1,F_2)&:M_1(i,\ell)=\ro^i \kl  i=\ell \kr,\\
	\Xi_4(F_1,F_2)&:\Fi(w,\tet,\ro,i)= \tet \ro ^i\lb 
\sum_{n=0}^{i-2} \ro ^{i(i-1)/2-n(n+1)/2} w^{i-1-n} \tet^n Q_n(w,\ro )Q(\ro ^{i}w,\ro )\right.\\
&\left.  + \sum_{n=i-1}^{\infty} \ro ^{n+1} w \tet^n Q_{i-2}(w,\ro )Q(\ro ^{n+2}w,\ro ) \rb,\mbox{ again, we have divided by }w,\\
\Xi_5(F_1,F_2)&=	\left\{ \sum_{u=2}^\II \pi(u)F_1^m(u,u-1)+\sum_{u=1}^\II \pi(u)\sum_{j=u}^\II   F_1^p(u,j)\right\}
\cdot \\ 
&\cdot \left\{ \sum_{\ell=2}^\II [M_1(j,\ell)+M_2(j,\ell)] C_2(l-1)  F_2^m(\ell,\ell -1)
+\sum_{\ell=1}^\II [M_1(j,\ell)+M_2(j,\ell)]\sum_{k=\ell}^\II  C_2(k)  F_2^p(\ell,k)\right\}.\\
&\mbox{This leads to }\\
\sig_Q^2&=\sig_3^2+2 \Xi_5(E_{w},E_{w})=0.4485678619\ldots,\\
\sig_X^2&=\sig_x^2+2\Xi_5(j,k)\equiv \sig_2^2,\\
C_{X,Q}&=\Eb(j\cdot E_{w}-\mu_2\mu_3)+\lb \Xi_5(j,E_{w})+\Xi_5(E_{w},k) \rb,\\
 \mu_3 &= 0.6188628379\ldots,\ro_{X,Q}=0.6289527540\ldots,\mu_4=1.6103718403\ldots,\sig_4^2=1.0188734817\ldots
\end{align*}
Again, we have made extensive simulations to check our results. The fit is quite good.
\section{The bargraph perimeter}
The wall (or bargraph) polyomino (wa) is made of contiguous columns, of any positive size,  with all base cells at the same level.

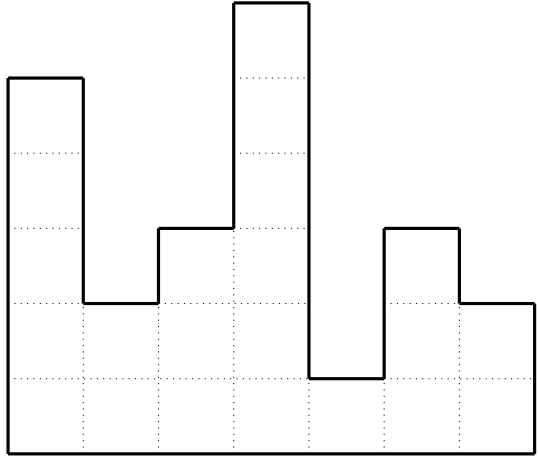
\begin{figure}[htbp]
\begin{center}
\begin{tikzpicture}[set style={{help lines}+=[dashed]}]
 \draw[very thick] (0,0) --  (7,0);
 \draw[very thick]  (0,0) --  (0,5);
 \draw[very thick]  (7,0) --  (7,2);
 \draw[very thick] (0,5) --  (1,5); \draw[very thick] (1,5) --  (1,2);
 \draw[very thick] (1,2) --  (2,2); \draw[very thick] (2,2) --  (2,3);
 \draw[very thick] (2,3) --  (3,3); \draw[very thick] (3,3) --  (3,6);
 \draw[very thick] (3,6) --  (4,6);
 \draw[very thick] (4,6) --  (4,1); \draw[very thick] (4,1) --  (5,1);
 \draw[very thick] (5,1) --  (5,3); \draw[very thick] (5,3) --  (6,3);
 \draw[very thick] (6,3) --  (6,2); \draw[very thick] (6,2) --  (7,2);
  \draw[dotted] (0,1) --  (4,1);  \draw[dotted] (5,1) --  (7,1);
  \draw[dotted] (0,2) --  (1,2);  \draw[dotted] (2,2) --  (4,2);  \draw[dotted] (5,2) --  (6,2);
  \draw[dotted] (0,3) --  (1,3);  \draw[dotted] (3,3) --  (4,3); 
  \draw[dotted] (0,4) --  (1,4);  \draw[dotted] (3,4) --  (4,4);  \draw[dotted] (3,5) --  (4,5); 
  \draw[dotted] (1,0) --  (1,2);  \draw[dotted] (2,0) --  (2,2); \draw[dotted] (3,0) --  (3,3);
  \draw[dotted] (4,0) --  (4,1);  \draw[dotted] (5,0) --  (5,1); \draw[dotted] (6,0) --  (6,2);
 \end{tikzpicture}
\end{center}
\caption{A typical wa  polyomino .}
\label{F10}
\end{figure}

 For wa, we have the following 
typical example: see Fig.\ref{F10},  

This polyomino is obviously equivalent to a composition of an integer $n$. In \cite{HL00}, it was proved that it asymptotically 
corresponds to a sequence of iid Geometric(1/2) RV.  We prove it again here, with a different method, more in the spirit of the present paper. 

\subsection{The generating functions }
We have
\bal
\FI(w,z)&=\sum_1^\II w^m\lp   \frac{z}{1-z}\rp^m=\frac{wz}{1-z-wz},\non \\
h(w,z)&=1-z-wz,\ro=\frac12,\non \\
[w^n z^n]&\FI(w,z)=\frac{1}{\ro^n}[w^n z^n]\FI(w,\ro z)=\frac{1}{\ro^n}[w^n z^n]\frac{wz\ro}{1-\ro z-\ro wz}.  \label{E5}
\end{align}
 We turn again to the analysis used by Bender in \cite{BE73}. It is asymptotically based on the GF
\[\FI_1(z)=\lp \frac{t(1)}{t(z)} \rp^m,\]
where $t(z)$ is the root of the denominator of (\ref{E5}), seen as a $w$ equation:
\[w=t(z)=\frac{1-\ro z}{\ro z},\FI_1(z)=\lp \frac{z/2}{1-z/2}\rp^m,\]
This indeed leads to sequence of iid Geometric(1/2) RV.
\subsection{The perimeter conditioned on $m$}

\[ \mbox{ with }x_d=j,x_{d-1}=k,\mbox{ we have here }w_d=|j-k|.\]
In \cite{GL18}, we have analyzed in great detail the perimeter of sequence of iid Geometric(p) RV, called therein a ``geometric word"
Using these results, we obtain
\bals
\mu_2&=\frac1p=2,\sig_x^2\equiv\sig_X^2\equiv\sig_2^2=\frac{1-p}{p^2}=2,\mu_3 =\frac{2(1-p)}{p(2-p)}=\frac43,
\sig_3^2= \frac{2(1-p)(p^2-2p+2)}{p^2(2-p)^2}=\frac{20}{9},\\
\sig_Q^2&=\frac{ 4(1-p)(p^4+9p^2-4p^3-10p+5)}{p^2(2-p)^2(p^2+3-3p)}=\frac{232}{63},\\
C_{X,Q}&=S_1+S_2-2\mu_2\mu_3,\\
S_1&=\sum_{i=1}^\II pq^{i-1}\lb \sum_{j=i}^\II pq^{j-1}j(j-i)+ \sum_{j=1}^{i-1} pq^{j-1}j(i-j) \rb,\\
S_2&=\sum_{i=1}^\II pq^{i-1}\lb \sum_{j=i}^\II pq^{j-1}i(j-i)+ \sum_{j=1}^{i-1} pq^{j-1}i(i-j) \rb,\\
C_{X,Q}&=\frac{ 2(2-4p+3p^2-p^3)}{p^2(2-p)^2}=\frac{20}{9},\\
\ro_{X,Q}&=\frac{5 \sqrt{2}\sqrt{406}}{174},\mu_4=5/3,\sig_4^2=\frac{173}{189}=.9153439162\ldots\\
\end{align*}
Again, we have made extensive simulations to check our results. The fit is quite good.
\subsection{A comparison with known GF}
Now we use  Bousquet-M\'elou \cite[(12)]{BM08}. If we set $q=z,y=1,x=w$ 
in the denominator $1-I_-$, we obtain a function
\[1-z-wz\]
which is exactly $h(w,z)$.

 If we set $q=z,y=x=v$,    we obtain
\[F(v,z)=1-\sum_{n=1}^\II \frac{ v^n (v-1)^j z^{j(j+1)/2}}{(z;z)_n (vz;z)_{n-1}}\]
which leads exactly, by Bender's theorems to $\mu_4,\sig_4^2$.
\section*{Acknowledgements}
We would like to thank  H.Prodinger and S.Wagner for solving a delicate recurrence.

\bibliographystyle{plain} 

\appendix
\section {Complements}
\subsection {$\Fi$-mixing property}                           \label{A1}
Let
\[\ldots,\xi_{-1},\xi_0,\xi_1,\ldots\] 
be a stationary sequence of RV. For $a\le b$, define $M_a^b$ as the $\sig$-field generated by $\xi_a,\ldots,\xi_b$,  define 
 $M_{-\infty}^a$ as the $\sig$-field generated by $\ldots,\xi_{a-1},\xi_a$, and define $M_a^\infty$ as the $\sig$-field generated by $\xi_a,\xi_{a+1},\ldots$.

Consider a nonnegative function $\Fi$ on positive integers, $\Fi(n)\ra 0,n\ra\infty$. We shall say that the sequence ${\xi_n}$ is $\Fi$-mixing if, for each $k(-\infty<k<\infty))$ and for each $n(n\ge 1), E_1\in M_{-\infty}^k,P(E_1)>0$, and $E_2 \in M_{k+n}^\infty$ together imply 
\[|\P(E_2|E_1)-\P(E_2)|\leq \Fi(n).\]
\subsection {Equation Bousquet-M\'elou \cite[(10)]{BM96}}  \label{A2}
Let $P$ be a polyomino. Let us mark by $s$ (resp.$t$) its left height (resp. right height), by  $x$ (resp. $y$)  the half-number of horizontal (resp. vertical) steps in its perimeter, by $q$ its area. The  dcc generating function $G(s,t,x,y,q)$ satisfies ($(a)_n$ means here $(a:q)_n$)
\bals
G(1,t,x,y,q)&=ty\frac{L_1(1)}{L_0(1)}, \mbox{ with }\\
L_0(s)&=1-\sum_{n=1}^\infty\frac{x^n s^n (y-1)^{n-1} q^{\bin{n+1}{2}}}   {(sq)_n (syq)_{n-1 } (syq)_n}, \mbox{ and }\\
L_1(s)&=\sum_{n=1}^\infty\frac{x^n s^n (y-1)^{n-1} q^{\bin{n+1}{2}}}   {(sq)_{n-1} (syq)^2_{n-1 } (1-styq^n)} ,\\
G(s,t,x,y,q)&=ty\frac{L_1(s)L_0(1)-L_1(1)L_0(s)+L_1(1)}{L_0(1)}.
\end{align*}
\subsection {Equation Bousquet-M\'elou \cite[(29)]{BM96}}  \label{A3}
The  cc generating function $G(s,t,x,y,q)$ satisfies
\bals
G(1,1,x,y,q)&=y\frac{(1-y)X}{1+W+yX},\\
X&=\frac{xq}{(1-y)(1-yq)}+\sum_{n=2}^\infty\frac{(-1)^{n+1} x^n (1-y)^{2n-4} q^{\bin{n+1}{2}} (y^2q)_{2n-2}     }  
 { (q)_{n-1} (yq)_{n-2} (yq)^2_{n-1}(yq)_n  (y^2q)_{n-1}   },\\
W&=  \sum_{n=1}^\infty  \frac{(-1)^{n} x^n (1-y)^{2n-3} q^{\bin{n+1}{2}} (y^2q)_{2n-1}        }  
 {  (q)_n(yq)^3_{n-1} (yq)_n (y^2q)_{n-1}     }.
\end{align*}
\subsection {Bender's theorems $1$ and $3$}  \label{A4}
Let $a_n(k)$ be a sequence of non-negative numbers. Let
\[p_n(k):=\frac{a_n(k)}{\sum_j a_n(j)}.\]
We say that $a_n(k)$ is asymptotically Gaussian (convergence in distribution)  with mean $\mu_n$ and variance $\sig^2_n$ if
\beq
\lim_{n\ra \infty}\sup_x \left| \sum_{k\leq \sig_n x+\mu_n} p_n(k)-\frac{1}{\sqrt{2\pi}}\int_\infty ^x e^{-t^2/2}dt \right|=0 .
\label{Ae1}
\eeq
We say that $a_n(k)$ satisfies a local limit theorem on a  set $S$ of real numbers if
\beq
\lim_{n\ra \infty}\sup_{x \in S}\left| \sig_n p_n(\sig_n x+\mu_n)-\frac{1}{\sqrt{2\pi}} e^{-x^2/2} \right|=0 .\label{Ae2}
\eeq
Theorem $1$: Central limit theorem

Let $f(z,w)$ have a power series expansion
\[f(z,w)=\sum_{n,k\geq 0}a_n(k)z^n w^k\]
with non-negative coefficients. Suppose there exists (i) an $A(s)$ continuous and non-zero near $0$, (ii) an $r(s)$ with bounded third derivative near $0$, (iii) a  non-negative integer $m$, and (iv) $\eps,\de>0$ such that
\[\lp 1-\frac{z}{r(s)}\rp^m f(z,e^s)-\frac{A(s)}{1-z/r(s)} \]
is analytic and bounded for
\[|s|<\eps,|z|<|r(0)|+\de.\]
Define 
\[\mu:=-\frac{r'(0)}{r(0)}, \sigd:=\mu^2-\frac{r''(0)}{r(0)}.\]
If $\sig \neq 0$, then (\ref{Ae1}) holds with $\mu_n=n\mu$ and $\sigd_n=n\sigd$.

\vspace{1ex}

Theorem $3$: Local limit theorem

Let $f(z,w)$ have a power series expansion
\[f(z,w)=\sum_{n,k\geq 0}a_n(k)z^n w^k\]
with non-negative coefficients and let $a<b$ be real numbers. Define
\[R(\eps):=\{z:a\leq \mbox{Re } z\leq b,|\mbox{Im } z|\leq \eps\}.\]
Suppose there exists $\eps>0,\de>0$,  a  non-negative integer $m$, and functions $A(s),r(s)$ such that
\bals
(i)& \hspace{1ex} A(s) \mbox{ is continuous and non-zero for } s\in R(\eps),\\
(ii)& \hspace{1ex}  r(s) \mbox{ is non-zero and has  a bounded third derivative for  } s\in R(\eps),\\
(ii)& \hspace{1ex}  \mbox{ for } s\in R(\eps) \mbox{ and } |z|<|r(s)|(1+\de)\\
& \hspace{1ex} \lp 1-\frac{z}{r(s)}\rp^m f(z,e^s)-\frac{A(s)}{1-z/r(s)}\\
&\mbox{ is analytic and bounded },\\
(iv)& \hspace{1ex} \lp \frac{r'(\al)}{r(\al)} \rp^2-  \frac{r''(\al)}{r(\al)} \neq 0 \mbox{ for  }a\leq \al \leq b,\\
(v)& \hspace{1ex} f(z,e^s) \mbox{ is analytic and bounded for}\\
&|z|\leq |r(\mbox{Re } s)|(1+\de) \mbox{ and }  \eps\leq |\mbox{Im } s| \leq \pi.
\end{align*}
Then we have
\[a_n(k)\sim \frac{n^m e^{-\al k}A(\al)}{m!r(\al)^n\sig_\al\sqrt{2\pi n}}\]
uniformly for $a\leq \al \leq b$, where
\bals
\frac{k}{n}&=-\frac{r'(\al)}{r(\al)},\\
\sigd_\al&=\lp \frac{k}{n} \rp^2- \frac{r''(\al)}{r(\al)}.
\end{align*}
In the proof of Theorem $3$, we find, on \cite[p.103]{BE73}), the following result: let $x:=(k-\mu_n)/\sig_n$, then
\[  \lim_{n\ra\infty}\left|\sig_n p_n(k)  -\frac{e^{-x^2/2}}{\sqrt{2\pi}}  \right| =0 \]
which is exactly the Local limit theorem (\ref{Ae2}) we need.

\end{document}